\documentclass [11pt,twoside,a4paper]{article}
\usepackage{amsfonts,epigraph}
\usepackage{amsthm}
\usepackage{amsmath,amssymb}
\usepackage{amstext}
\usepackage{url}
\usepackage{nicefrac}
\usepackage{mathrsfs}
\usepackage{amscd}
\usepackage{xypic}
\usepackage{epsf} %图形宏包
\usepackage{graphicx} %这是图形命令的宏包
\usepackage{fancybox} %这是盒子命令的宏包
\usepackage{color} %这是颜色宏包
\usepackage{fancyhdr}
\usepackage[hang,footnotesize]{caption2} %有关图形标题的宏包
\usepackage[outercaption]{sidecap}
\usepackage{booktabs} %for table and tabular
\usepackage[all,cmtip]{xy}%复杂交换图

\newfont{\aaa}{cmb10 at 19pt}
\newfont{\bbb}{cmb10 at 11pt}
\newtheorem{theorem}{Theorem}[section]
\newtheorem{corollary}[theorem]{Corollary}
\newtheorem{lemma}[theorem]{Lemma}
\newtheorem{proposition}[theorem]{Proposition}

\theoremstyle{definition}
\newtheorem{definition}[theorem]{Definition}

\newtheorem{ass}[theorem]{Assumption}

\newtheorem{choices-notations}[theorem]{Choices and Notations}

\newtheorem{rem}[theorem]{Remark}

\def\v1{\vspace{1mm}}

\def\leq{\leqslant}

\def\geq{\geqslant}
\def\NN{\Bbb N}
\def\RR{\mathbb R}
\def\ZZ{\mathbb Z}
\def\CC{{\mathbb C}}

\theoremstyle{remark}

\newcommand{\exendproof}{\renewcommand{\qed}{\relax}\end{proof}}

\newsavebox{\SmallMathBox}

\DeclareRobustCommand*{\nicefrac}[2]{\ifmmode\mathnicefrac{#1}
{ #2}%
  \else\textnicefrac{#1}{#2}\fi}
\newcommand*{\textnicefrac}[2]{\check@mathfonts%
\mbox{\raisebox{.5ex}{\fontsize\sf@size\z@\selectfont#1}\kern-.
1em%
/\kern-.1em\raisebox{- .25ex}{\fontsize\sf@size\z@\selectfont#2} }}
\newcommand*{\mathnicefrac}[2]{%
  \mathchoice
    {\m@fr@c{\scriptstyle}{#1}{#2}}
    {\m@fr@c{\scriptstyle}{#1}{#2}}
    {\m@fr@c{\scriptscriptstyle}{#1}{#2}}
    {\m@fr@c{\scriptscriptstyle}{#1}{#2}}}
\def\sqm1{\sqrt{-1}}

\def\={\cong}
\def\>{\supset}
\def\<{\subset}

\def\12{\frac{1}{2}}
\def\0{^{\circ}}

\def\Bb{{\mathcal B}}

\def\Ll{{\mathcal L}}

\def\Pp{{\mathcal P}}

\DeclareMathOperator{\diag}{diag}
\DeclareMathOperator{\End}{End}

\DeclareMathOperator{\GL}{GL}  \DeclareMathOperator{\Graph}{Gr}
 
\DeclareMathOperator{\Hom}{Hom}

\DeclareMathOperator{\Mas}{Mas}

\DeclareMathOperator{\sign}{sign}

 \DeclareMathOperator{\Sp}{Sp}

\DeclareMathOperator{\ssp}{sp}

\makeatletter
\def\@evenhead{
\vbox{\hbox to \textwidth {}{\hspace{0mm}{\footnotesize
\thepage}}{\hspace{9.8cm} {\footnotesize {Yuting ZHOU et al.}}} \protect\vspace{1truemm}\relax \hrule depth0pt
height0.15truemm width\textwidth}}
\def\@evenfoot{}
\def\@oddhead{\vbox{\hbox to \textwidth
{{\hspace{0cm}{\footnotesize The H{\"o}rmander index in the finite-dimensional case}\hfill{\footnotesize
\thepage}}\hspace{0mm}}{} \protect\vspace{1truemm}\relax\hrule
depth0pt height0.15truemm width\textwidth}}
\def\@oddfoot{}
\makeatother
\begin{document}
%首页眉设定%%%%%%%%%%%%%%%%%%%%%%%%%%%%%%%%%%%%%%%%%%%%%%%%%%
%\setcounter{page}{61} %
\thispagestyle{empty} \thispagestyle{fancy} {
\fancyhead[lO,RE]{\footnotesize Front. Math. China \\
DOI 10.1007/s11464-017-0686-5\\[3mm]
}
\fancyhead[RO,LE]{\scriptsize \bf %http://www.spcum.hep.com.cn\\
} \fancyfoot[CE,CO]{}}
\renewcommand{\headrulewidth}{0pt}%(首页书眉线宽）
%\renewcommand{\headsep}{0.7cm}%（首页正文与书眉距离）

%首页眉设定%%%%%%%%%%%%%%%%%%%%%%%%%%%%%%%%%%%%%%%%%%%%%%%%%%

%%%以下正文开始
\setcounter{page}{1}
\qquad\\[8mm]

\noindent{\aaa{The H{\"o}rmander index in the finite-dimensional case}}\\[1mm]

\noindent{\bbb Yuting ZHOU$^1$,\quad Li WU$^2$,\quad Chaofeng ZHU$^1$}\\[-1mm]

\noindent\footnotesize{1 Chern Institute of Mathematics and LPMC,
Nankai University, Tianjin 300071, China}

\noindent\footnotesize{2 Department of  Mathematics,
Shandong University, Jinan, 250100, China}\\[6mm]

\vskip-2mm \noindent{\footnotesize$\copyright$ Higher Education
Press and Springer-Verlag Berlin Heidelberg 2017} \vskip 4mm

\normalsize\noindent{\bbb Abstract}\quad In this paper, we calculate H{\"o}rmander index in the finite-dimensional case. Then we use the result to give some iteration
inequalities, and prove almost existence of mean indices for given
complete autonomous Hamiltonian system on compact symplectic
manifold with symplectic trivial tangent bundle and given autonomous
Hamiltonian system on regular compact energy hypersurface of
symplectic manifold with symplectic trivial tangent bundle.\vspace{0.3cm}

\footnotetext{Received December 14, 2016; accepted January 18,
2017\\
\hspace*{5.8mm}Corresponding author: Chaofeng ZHU, E-mail:
zhucf@nankai.edu.cn}

\noindent{\bbb Keywords}\quad Maslov index, H{\"o}rmander index, Maslov-type index,
symplectic reduction\\
{\bbb MSC}\quad 53D12, 58J30\\[0.4cm]
%\noindent{\bbb 1\quadIntroduction}
\section{Introduction}\label{s:introduction}

Let $(V,\omega)$ be a symplectic vector space. Let $\lambda_{1},
\lambda_{2}, \mu_{1}, \mu_{2}$ be four Lagrangian subspaces of $V$.
The H{\"o}rmander index $s(\lambda_{1},\lambda_{2};\mu_{1}, \mu_{2}
)$ has been introduced by L. H{\"o}rmander \cite[Sect. 3.3]{Ho71}
in the finite-dimensional case when $H = \lambda_{j}\oplus \mu_{k}$
hold for $j,k =1,2$, who also gave the explicit formula to calculate
it. The notion was generalized by B. Booss and K. Furutani
\cite[Proposition 2.1]{BoFu99} in the finite-dimensional case and
\cite[Definition5.2]{BoFu99} in the strong symplectic Hilbert case
when $(\lambda_{1},\mu_{1}), (\lambda_{2},\mu_{1})$ are two Fredholm
pairs of Lagrangian subspaces of V and
$\mu_{1}/(\mu_{1}\cap\mu_{2})$ is finite-dimensional. Recently, B.
Booss and the third author \cite[Definition 3.4.4]{BoZh14}
generalized the notion to the symplectic Banach case.

The splitting number \cite[Denition 9.1.4]{Lo02} is a special case
of the H{\"o}rmander index. It turns out that the study of the
H{\"o}rmander index in the full generality is very important in the
study of Hamiltonian systems (see \cite{Lo02} for the applications
of the splitting numbers and \cite{LoZhZh06, LiZh14} for the study
of multiplicity of the brake orbits).

In \cite{Go09}, M. de Gosson gave a very elegant definition of the
H{\"o}rmander index in the finite-dimensional case in great
generality. His definition differs slightly from ours. By admitting
half-integer indices, it yields more simple proofs, but may be more
difficult to be used in concrete applications in Morse theory.

We calculate the H{\"o}rmander index in the finite-dimensional case
and get the following main result.

\begin{theorem}\label{t:calculate-hormander-index}
Let $(V,\omega)$ be a complex symplectic vector space  of dimension
$2n$. Let $\lambda_{1}, \lambda_{2}, \mu_{1}, \mu_{2}$ be four
Lagrangian subspaces of $V$. Denote by $i(\cdot,\cdot,\cdot)$ the
triple index defined by \cite[(2.16)]{Du76} (see Corollary
\ref{c:i-index-definition} below). Then we have
\begin{eqnarray}
% \nonumber to remove numbering (before each equation)
  s(\lambda_{1},\lambda_{2};\mu_{1},\mu_{2}) &=& i(\lambda_{1},\lambda_{2},\mu_{2})-i(\lambda_{1},\lambda_{2},\mu_{1}) \\
   &=&i(\lambda_{1},\mu_{1},\mu_{2})-i(\lambda_{2},\mu_{1},\mu_{2}).
\end{eqnarray}
\end{theorem}

Our main result does not require the transversal conditions $V
=\lambda_{j}\oplus\mu_{k}, j,k = 1,2$. We use the result to get some
new iteration inequalities of Maslov-type index. Then we use the
inequalities to prove almost existence of mean indices for given
complete autonomous Hamiltonian system on compact symplectic
manifold with symplectic trivial tangent bundle and given autonomous
Hamiltonian system on regular compact energy hypersurface of
symplectic manifold with symplectic trivial tangent bundle.

The paper is organized as follows. In $\S \ref{s:introduction}$, we
review the historical literatures and introduce our main result. In
$\S \ref{s:Maslov-index}$, we review the notions of Maslov index and
Maslov-type index in the finite-dimensional case. In $\S
\ref{s:calculation}$, we study the H{\"o}rmander index in the finite-dimensional case and prove Theorem
\ref{t:calculate-hormander-index}. In $\S
\ref{s:iteration-inequalities}$, we prove some iteration
inequalities of Maslov-type index. In $\S \ref{s:mean-indices}$, we
prove almost existence of mean indices for given complete autonomous
Hamiltonian system on compact symplectic manifold with symplectic
trivial tangent bundle and given autonomous Hamiltonian system on
regular compact energy hypersurface of symplectic manifold with
symplectic trivial tangent bundle.

In this paper we denote the sets of natural, integral, real, complex
numbers, the unit circle in the complex plane, the set of all linear
operators on V, the set of all invertible linear transformations on
$V$ and the set of all self-adjoint operators on Hilbert space $V$,
by $\NN$, $\ZZ$, $\RR$, $\CC$, $S^{1}$, $\End(V)$, $\GL(V)$ and
$\Bb^{sa}(V)$ respectively. We denote by $\dim V$ the
complex dimension of a complex linear space $V$. We denote by
$\Hom(V_{1},V_{2})$ the set of all linear maps between vector spaces
$V_{1}$ and $V_{2}$. For  a linear operator $A\in
\Hom(V_{1},V_{2})$, we denote by $\Graph(A)$ the graph of $A$. For
two maps $f: X \rightarrow Y$ and $g: Y\rightarrow Z$, we denote by
$g\circ f: X\rightarrow Z$ the composite map defined by $(g\circ
f)(x)=g(f(x))$ for each $x\in X$. Without further explanation, the
coefficient field is $\CC$ in the rest of this paper.

%\noindent\\[4mm]

%\noindent{\bbb 2\quad The Maslov index and the Maslov-type index}\\[0.1cm]
\section{The Maslov index and the Maslov-type index}\label{s:Maslov-index}
In this section we review the definition of the Maslov index and the
Maslov-type index in the finite-dimensional case. Firstly we recall
the basic concepts and properties of symplectic vector space.

\begin{definition}\label{d:symplectic-space}
Let $V$ be a complex vector space.
\begin{enumerate}
\item[(a)]  A mapping
\begin{equation*}
    \omega: V\times V\rightarrow \CC
\end{equation*}
is called a \textit{symplectic form} on $V$, if it is a
non-degenerate skew-Hermitian form. Then we call $(V,\omega)$ a
\textit{complex symplectic vector space}.

\item[(b)]  Let $(V,\omega_1)$ and $(V,\omega_2)$ be two finite-dimensional symplectic vector spaces. A linear map $L\in
\Hom(V_{1},V_2)$ is called \textit{symplectic}, if $L$ is invertible
and $\omega_2(Lx,Ly)=\omega_1(x,y)$ for each $x,y\in V_1$. We denote
by $\Sp((V_1,\omega_1),(V_2,\omega_2))$ the set of all such
symplectic linear maps $L$, and
$\Sp(V,\omega):=\Sp((V,\omega),(V,\omega))$ for the finite-dimensional symplectic vector space $(V,\omega)$. We denote by
$\Sp(V_1,V_2)=\Sp((V_1,\omega_1),(V_2,\omega_2))$ and
$\Sp(V)=\Sp(V,\omega)$ if there is no confusion.

\item[(c)]  Let $\lambda$ be a linear subspace of $V$. The
\textit{annihilator} $\lambda^{\omega}$ of $\lambda$ is defined by
\begin{equation*}
\lambda^{\omega}:=\{x\in V; \omega(x,y)=0\; \text{for all} \; y\in
\lambda\}.
\end{equation*}

\item[(d)]  A linear subspace $\lambda$ of $V$ is called \textit{symplectic}, \textit{isotropic}, \textit{co-isotropic},
or \textit{Lagrangian} if
\begin{equation*}
    \lambda\cap\lambda^{\omega}=\{0\},\;\;\;\lambda\subset\lambda^{\omega},\;\;\;\lambda\supset\lambda^{\omega},\;\;\;\lambda=\lambda^{\omega},
\end{equation*}
respectively.

\item[(e)]  The \textit{Lagrangian Grassmannian}
$\Ll(V,\omega)$ consists of all Lagrangian subspaces of
$(V,\omega)$. We write $\mathcal{L}(V):=\mathcal{L}(V,\omega)$ if
there is no confusion.
\end{enumerate}
\end{definition}
If $\dim V < +\infty$, the space $\mathcal{L}(V)$ is path-connected.
It is nonempty if and only if the signature $\sign(i\omega)=0$.

The following definition of the Maslov index is taken from
\cite[\S2.2]{BoZh14}.

Let $(V,\omega)$ be a $2n$-dimensional complex vector space. Let
$\langle\cdot,\cdot\rangle$ be an inner product on $V$. Then there
is an operator $J\in \GL(V)$ such that $iJ$ is self-adjoint and
$\omega(x,y)=\langle Jx,y\rangle$ for all $x,y \in V$.

Denote by $V^{\mp}$ the positive (negative) eigenspace of $iJ$.
Given a path $(\lambda(s),\mu(s))$, $s\in[a,b]$ of pairs of
Lagrangian subspaces of $(V,\omega)$, let $U(s), V(s) :
V^-\rightarrow V^+$ be generators for $(\lambda(s),\mu(s))$,
i.e., $\lambda(s)=\Graph(U(s))$ and $\mu(s)=\Graph(V(s))$ (see
\cite[Proposition 2]{BoZh13}). Then the family
$\{U(s)V(s)^{-1}\}_{s\in[a,b]}$ is a continuous family of unitary
operators on Hilbert space $(V^+,-i\omega|_{V^+})$.

Note that the eigenvalues of $U(s)V(s)^{-1}$ are on the unit circle
$S^1$. Recall that each map in $C([a,b], S^1)$ can be lifted to a
map in $C([a,b],\RR)$. By \cite[Theorem II.5.2]{Ka95}, there are $n$
continuous functions $\theta_1,\ldots, \theta_n \in C([a,b],\RR)$
such that the eigenvalues of the operator $U(s)V(s)^{-1}$ for each
$s\in [a,b]$ (counting algebraic multiplicity) have the form
\begin{equation*}
    e^{i\theta_j(s)},\; j=1,\ldots,n.
\end{equation*}

Denote by $[a]$ the integer part of $a\in\RR$ and $\{a\}:= a-[a]$.
Define
\begin{equation}\label{e:integerE}
E(a):=
\begin{cases}
a,&a\in\ZZ;\\
[a]+1,&a\not\in\ZZ.
\end{cases}
\end{equation}

\begin{definition}\label{d:maslov-index}
We define the \textit{Maslov} \textit{index} of the path
$(\lambda,\mu)$ by
\begin{equation}\label{e:maslov+}
    \Mas\{\lambda,\mu\}=\Mas_{+}\{\lambda,\mu\}=\sum_{j=1}^{n}\biggl( E(\frac{\theta_j(b)}{2\pi})-E(\frac{\theta_j(a)}{2\pi})\biggr),
\end{equation}

\begin{equation}\label{e:maslov-}
\Mas_-\{\lambda,\mu\}=\sum_{j=1}^n\biggl([\frac{\theta_j(b)}{2\pi}]-[\frac{\theta_j(a)}{2\pi}]\biggr).
\end{equation}
\end{definition}

By definition, $\Mas_{\pm}\{\lambda(s),\mu(s);s\in[a,b]\}$ is an
integer which does not depend on the choices of the arguments
$\theta_j(s)$. By \cite[Proposition 6]{BoZh13}, it does not depend
on the particular choice of the inner product.

\begin{rem}\label{r:maslov-index}
Let $V$ be a real symplectic $2n$-dimensional symplectic space. Let
$\alpha$ be in $\mathcal{L}(V)$ and $\lambda(s)$, $s\in[a,b]$ be a
path in $\mathcal{L}(V)$ such that
$\lambda(a)\cap\alpha=\lambda(b)\cap\alpha=\{0\}$. Denote by
$[\lambda:\alpha]$ the Maslov-Arnol'd index defined by
\cite[(2.8)]{Du76}. By \cite[Proposition 3.27]{BoZh14}, we have
$\Mas\{\lambda,\alpha\}=-\Mas\{\alpha,\lambda\}=-[\lambda:\alpha].$
\end{rem}

We recall some concepts in \cite{Du76,Zh06} for the calculation of
the Maslov index.

Let $(V,\omega)$ be a finite-dimensional symplectic vector space.
Let $\alpha(s)$, $s\in(-\varepsilon,\varepsilon)$ be a path in
$\mathcal{L}(V)$ differentiable at $s=0$. Define the form
$Q(\alpha,0)$ on $\alpha(0)$ by

\begin{equation}\label{e:Q-form}
Q(\alpha,0)(x,y):=\frac{d}{ds}|_{s=0}\omega(x,y(s)),
\end{equation}
where $x,y\in \alpha(0)$, $y(s)\in \alpha(s)$ and $y(s)-y \in
\beta$. It is well-known that the Hermitian form $Q(\alpha,0)$ is
independent on the choice of $\beta\in\mathcal{L}(V)$ with
$V=\alpha(0)\oplus\beta$.

\begin{definition}\label{d:crossing-form}
Let $(V,\omega)$ be a finite-dimensional symplectic vector space.
Let $(\lambda(s),\mu(s))$, $s\in [a,b]$ be a $C^1$ curve of pairs of
Lagrangian subspaces of $(V,\omega)$. For $t\in [a,b]$, the
\textit{crossing  form} $\Gamma(\lambda,\mu ,t)$ on
$\lambda(t)\cap\mu(t)$ is defined by
\begin{equation}\label{e:crossing-form}
\Gamma(\lambda,\mu,t)(u,v):=Q(\lambda,t)(u,v)-Q(\mu,t)(u,v),
\end{equation}
where $u,v \in \lambda(t)\cap \mu(t)$. A \textit{crossing} is a time
$t\in[a,b]$ such that $\lambda(t)\cap\mu(t)\neq\{0\}$. A crossing
$t$ is called \textit{regular} if $\Gamma(\lambda,\mu,t)$ is
non-degenerate.
\end{definition}

\begin{proposition}\label{p:calculate-maslov-index}(\cite[Lemma 2.5]{Du76}, \cite[Proposition 4.1]{Zh06})
Let $(V,\omega)$ be a finite-dimensional symplectic vector space.
Let $(\lambda(s),\mu(s))$, $s\in[a,b]$ be a $C^1$ curve of pairs of
Lagrangian subspaces of $(V,\omega)$ with only regular crossings.
Then the crossings are finite, and we have
\begin{equation}\label{e:calculate-maslov-index}
\Mas\{\lambda,\mu\}=m^{+}(\Gamma(\lambda,\mu,a))-m^{-}(\Gamma(\lambda,\mu,b))+\sum_{a<s<b}\sign(\Gamma(\lambda,\mu,s)),
\end{equation}
where we denote by $m^*(Q)$, $*=+,0,-$ the Morse positive index, the
nullity, and the Morse negative index of an Hermitian form $Q$
respectively.
\end{proposition}

For each $\tau >0$, we define
\begin{equation}\label{e:symplectic-path-I}
\mathcal{P}_{\tau}(V):=\{\gamma\in C([0,\tau],\Sp(V));
\gamma(0)=I_{V}\}.
\end{equation}

\begin{definition}\label{d:Maslov-type-index}(cf. \cite[Definition 4.6]{Zh06})
Let $(V_l,\omega_l)$, $l=1,2$ be two finite-dimensional symplectic
vector space. Then $(V=V_1\oplus V_2, (-\omega_1)\oplus\omega_2)$ is
a symplectic vector space. Let $W\in \mathcal{L}(V)$. Let
$\gamma(t)$, $0\leq t\leq\tau$ be a path in $\Sp(V_1,V_2)$. The
\textit{Maslov-type index} $i_W(\gamma)$ is defined to be
$\Mas\{\Graph\circ\gamma, W\}$. If $P\in \Sp(V_1,V_2)$, we define
$i_P(\gamma):=i_{\Graph(P)}(\gamma)$. If
$(V_1,\omega_1)=(V_2,\omega_2)$ and
$\gamma\in\mathcal{P}_{\tau}(V_1)$, we denote by
$i_1(\gamma):=i_{I_{V_1}}(\gamma)$, $\nu_1(\gamma):=\dim
\ker(\gamma(\tau)-I_{V_1})$, $\nu_{W}(\gamma):=\dim
(\Graph(\gamma(\tau))\cap W)$ and $\nu_{P}(\gamma):=\dim \ker
(\gamma(\tau)-P)$.
\end{definition}

The following lemma gives the tangent vector of given symplectic
path.

\begin{lemma}\label{l:tangential-map}(cf. \cite[Lemma 3.1]{Du76})
Let $(V_l,\omega_l)$, $l=1,2$ be two finite-dimensional symplectic
vector space. Then $(V=V_1\oplus V_2,(-\omega_1)\oplus \omega_2)$ is
a symplectic vector space. Let $\gamma(t)$, $t\in(-\varepsilon,
\varepsilon)$ be a path in $\Sp(V_1,V_2)$. Then for each $u,v\in
V_1$, we have
\begin{equation}\label{e:tangential-map}
Q(\Graph
\circ\gamma,0)((u,\gamma(0)u),(v,\gamma(0)v))=\omega_1(-\gamma(0)^{-1}\dot{\gamma}(0)u,v),
\end{equation}
where we denote by $\dot{\gamma}:=\frac{d}{dt}\gamma$.
\end{lemma}

The Maslov-type indices have the following properties.

\begin{lemma}\label{l:Maslov-index-neighborhood}(cf. \cite[Theorem 6.1.8]{Lo02})
Let $(V,\omega)$ be a finite-dimensional symplectic vector space.
Let $W\in \mathcal{L}(V\oplus V,(-\omega)\oplus \omega)$. Let
$\gamma(t)\in\mathcal{P}_{\tau}(V)$ be a symplectic path. Then there
exists a $C^0$ neighborhood $\mathcal{N}$ of $\gamma$ in
$\mathcal{P}_{\tau}(V)$ such that, for each $\tilde{\gamma}\in
\mathcal{N}$, there holds that
\begin{equation}\label{e:Maslov-index-neighborhood}
i_W(\gamma)\leq i_W(\tilde{\gamma}) \leq i_W(\gamma)+\nu_W(\gamma).
\end{equation}
\end{lemma}

\begin{proof}
Let $N$ be a convex neighborhood of $0$ in $\ssp (V)$, the Lie
algebra of $\Sp(V)$ such that the map $\exp: N\rightarrow \exp(N)$
is a diffeomorphism. Set
\begin{equation}\label{e:neighborhood}
\mathcal{N}:=\{\tilde{\gamma}\in \mathcal{P}_{\tau}(V);
\tilde{\gamma}(t)\in \gamma(t)\exp(N),\  \text{for each}\  t\in
[0,\tau] \}.\nonumber
\end{equation}
Then for each $\tilde{\gamma}\in \mathcal{N}$, there exists $A\in
C([0,\tau],N)$ and $A(0)=0$, such that
$\tilde{\gamma}(t)=\gamma(t)\exp(A(t))$, $t\in[0,\tau]$. Thus we
have a homotopy $\phi(s,t):[0,1]\times[0,\tau]\rightarrow \Sp(V)$ by
$\phi(s,t)=\gamma(t)\exp(sA(t))$. Denote by $\gamma_{\tau}$ the path
$\gamma(\tau)\exp(sA(\tau))$, $s\in [0,1]$. By the definition of the
Maslov-type index we can choose $N$ small enough such that $0\leq
i_W(\gamma_{\tau})\leq\nu_W(\gamma)$ for each $\tilde{\gamma}\in
\mathcal{N}$. By the homotopy invariance and path additivity of the
Maslov-type index we have
$i_W(\tilde{\gamma})=i_W(\gamma)+i_W(\gamma_{\tau})$. The inequality
(\ref{e:Maslov-index-neighborhood}) then follows.
\end{proof}

\begin{lemma}\label{l:calculate-Maslov-multiply}(cf. \cite[Lemma 4.4]{Zh06})
Let $(V_l,\omega_l)$, $l=1,2,3,4$ be finite-dimensional symplectic
vector spaces. Let $W$ be a Lagrangian subspace of $(V_1\oplus
V_4,(-\omega_1)\oplus \omega_4)$. Let $\gamma_l\in
C([0,1],\Sp(V_l,V_{l+1}))$, $l=1,2,3$ be symplectic paths. Then we
have
\begin{equation}\label{e:calculate-Maslov-multiply}
i_W(\gamma_3\gamma_2\gamma_1)=i_{W'}(\gamma_2)+i_W(\gamma_3\gamma_2(0)\gamma_1),
\end{equation}
where $W'=\diag (\gamma_1(1),\gamma_3(1)^{-1})W$.
\end{lemma}

\section{Calculation of the H{\"o}rmander index}\label{s:calculation}
In this section we study the form $Q(\cdot,\cdot;\cdot)$ and the
triple index $i(\cdot,\cdot,\cdot)$. Then we express the
H{\"o}mander index via the triple index.

\subsection{The form $Q(\cdot,\cdot;\cdot)$}\label{s:form-Q}
Let $(V,\omega)$ be a complex symplectic vector space with three
isotropic subspaces $\alpha,\beta,\gamma$. Define
$Q:=Q(\alpha,\beta;\gamma)$ on $\alpha\cap(\beta+\gamma)$ by
\begin{equation}\label{e:Q-definition}
Q(x_1,x_2):=\omega(x_1,y_2)=\omega(z_1,y_2)=\omega(x_1,z_2)
\end{equation}
for all $x_j=-y_j+z_j\in\alpha\cap(\beta+\gamma)$, where
$y_j\in\beta$, $z_j\in\gamma$, $j=1,2$.

\begin{rem}\label{r:Q-Duistermaat}
Assume that $\alpha,\beta,\gamma$ are Lagrangian subspaces of $V$
and $V=\alpha\oplus\beta=\beta\oplus\gamma$. Then the form
$Q(\alpha,\beta;\gamma)$ defined by \cite[(2.3)]{Du76} is
$-Q(\alpha,\beta;\gamma)$ here.
\end{rem}

With the above notions, we have $\omega(x_1+y_1,x_2+y_2)=0$,
$\omega(x_1,x_2)=0$ and $\omega(y_1,y_2)=0$. It follows that
$\omega(x_1,y_2)=-\omega(y_1,x_2)=\overline{\omega(x_2,y_1)}$. So
$Q(x_1,x_2)$ does not depend on the choices of $y_1$ and $y_2$ and
the form $Q$ is a well-defined Hermitian form on
$\alpha\cap(\beta+\gamma)$. Moreover, the last two equalities in
(\ref{e:Q-definition}) hold, i.e.,
$\omega(x_1,y_2)=\omega(z_1,y_2)=\omega(x_1,z_2)$, and we have
\begin{equation}\label{e:Q-transposition}
Q(\alpha,\beta;\gamma)=-Q(\alpha,\gamma;\beta).
\end{equation}

The following lemma is well-known in the non-degenerate case.

\begin{lemma}\label{l:Q-cyclic}
Let $(V,\omega)$ be a complex symplectic vector space with three
isotropic subspaces $\alpha,\beta,\gamma$. Denote by
$Q_1:=Q(\alpha,\beta;\gamma)$, $Q_2:=Q(\beta,\gamma;\alpha)$ and
$Q_3=Q(\gamma,\alpha;\beta)$.
\newline (a) Let $z_j=x_j+y_j\in \gamma\cap(\alpha+\beta)$, where
$x_j\in\alpha$, $y_j\in \beta$, $j=1,2$. Then we have
$Q_1(x_1,x_2)=Q_2(y_1,y_2)=Q_3(z_1,z_2)$.
\newline (b) We have $m^{\pm}(Q_1)=m^{\pm}(Q_2)=m^{\pm}(Q_3)$.
\end{lemma}
\begin{proof}
(a) Note that $x_j=z_j-y_j$. So we have
\begin{equation*}
    Q_1(x_1,x_2)=\omega(x_1,z_2)=\omega(y_1,-z_2)=Q_2(y_1,y_2).
\end{equation*}
Similarly we have $Q_1(x_1,x_2)=Q_3(z_1,z_2)$.

(b) Let $X$ be a linear subspace of $\alpha\cap(\beta+\gamma)$ such
that $Q_1|_{X}<0$ with a base $\{x_1,...,x_k\}$. Let $y_j\in \beta$,
$j=1,...,k$ be such that $x_j+y_j\in\gamma$. Then $y_j\in
\beta\cap(\alpha+\gamma)$. Let $Y$ be the $\CC$-linear span of
$\{y_1,...,y_k\}$. By (a), for each $(a_1,...,a_k)\in
\CC^k\backslash\{0\}$ we have $Q_2(y,y)=Q_1(x,x)<0$, where
$x:=\sum_{j=1}^{k}a_jx_j$ and $y:=\sum_{j=1}^k a_jy_j$. In this case
$y\neq0$. So $y_1,...,y_k$ are linearly independent, and $\dim Y=k$.
Hence $m^-(Q_2)\geq \dim Y=\dim X$. Since $X$ is arbitrarily chosen,
we have $m^-(Q_2)\geq m^-(Q_1)$. Similarly we have $m^-(Q_1)\geq
m^-(Q_2)$. It follows that $m^-(Q_1)=m^-(Q_2)$. Similarly we have
$m^+(Q_1)=m^+(Q_2)$ and $m^{\pm}(Q_2)=m^{\pm}(Q_3)$.
\end{proof}

Here we give the kernel of the form $Q(\alpha,\beta;\gamma)$.
\begin{lemma}\label{l:kerQ-isotopic}
Let $(V,\omega)$ be a complex symplectic vector space with three
isotropic subspaces $\alpha,\beta,\gamma$. Then we have
\begin{eqnarray}
% \nonumber to remove numbering (before each equation)
  \alpha\cap(\beta\cap(\alpha+\gamma))^{\omega} &=& \alpha\cap(\gamma\cap(\alpha+\beta))^{\omega},\  and \label{e:kerQ-00}\\
  \ker Q(\alpha,\beta;\gamma) &=& \alpha\cap(\beta+\gamma\cap(\beta\cap(\alpha+\gamma))^{\omega}) \label{e:kerQ-1} \\
   &=& \alpha\cap(\beta+\gamma\cap(\alpha\cap(\beta+\gamma))^{\omega}) \label{e:kerQ-2} \\
   &=&\alpha\cap(\beta\cap(\gamma\cap(\alpha+\beta))^{\omega}+\gamma).
   \label{e:kerQ-3}
\end{eqnarray}
In particular, $\ker Q(\alpha,\beta;\gamma)=\alpha\cap\beta
+\alpha\cap\gamma$ holds if
$(\alpha\cap(\beta+\gamma))^{\omega}=\alpha+\beta\cap\gamma$, or
$(\beta\cap(\alpha+\gamma))^{\omega}=\beta+\alpha\cap\gamma$, or
$(\gamma\cap(\alpha+\beta))^{\omega}=\gamma+\alpha\cap\beta$.
\end{lemma}
\begin{proof}
Since $\beta\cap(\alpha+\gamma)\subset
\alpha+\gamma\cap(\alpha+\beta)$, we have
\begin{eqnarray*}
% \nonumber to remove numbering (before each equation)
  \alpha\cap(\beta\cap(\alpha+\gamma))^{\omega} &\supset &\alpha\cap(\alpha+\gamma\cap(\alpha+\beta))^{\omega} \\
   &=&\alpha\cap(\alpha^{\omega}\cap(\gamma\cap(\alpha+\beta))^{\omega}) \\
   &=& \alpha\cap(\gamma\cap(\alpha+\beta))^{\omega}.
\end{eqnarray*}
Similarly we have
$\alpha\cap(\beta\cap(\alpha+\gamma))^{\omega}\subset
\alpha\cap(\gamma\cap(\alpha+\beta))^{\omega}$. So (\ref{e:kerQ-00}) holds.

Since
$\beta\subset\beta^{\omega}\subset(\beta\cap(\alpha+\gamma))^{\omega}$,
by (\ref{e:Q-definition}) we have
\begin{eqnarray*}
% \nonumber to remove numbering (before each equation)
  \ker Q(\alpha,\beta;\gamma) &=& \alpha\cap(\beta+\gamma)\cap(\beta\cap(\alpha+\gamma))^{\omega} \\
   &=&
   \alpha\cap(\beta+\gamma\cap(\beta\cap(\alpha+\gamma))^{\omega}).
\end{eqnarray*}
So (\ref{e:kerQ-1}) follows. Similarly we get (\ref{e:kerQ-3}). By
(\ref{e:kerQ-00}) and (\ref{e:kerQ-1}) we get (\ref{e:kerQ-2}).

If $(\alpha\cap(\beta+\gamma))^{\omega}=\alpha+\beta\cap\gamma$, by
(\ref{e:kerQ-2}) we have
\begin{eqnarray*}
% \nonumber to remove numbering (before each equation)
  \ker Q(\alpha,\beta;\gamma) &=& \alpha\cap(\beta+\gamma\cap(\alpha+\beta\cap\gamma)) \\
   &=& \alpha\cap(\beta+\alpha\cap\gamma +\beta\cap\gamma) \\
  &=& \alpha\cap(\beta+\alpha\cap\gamma)=\alpha\cap\beta +
  \alpha\cap\gamma.
\end{eqnarray*}
Similarly, by (\ref{e:kerQ-1}) and (\ref{e:kerQ-3}), $\ker
Q(\alpha,\beta;\gamma)=\alpha\cap\beta+\alpha\cap\gamma$ holds if
$(\beta\cap(\alpha+\gamma))^{\omega}=\beta+\alpha\cap\gamma$, or
$(\gamma\cap(\alpha+\beta))^{\omega}=\gamma+\alpha\cap\beta$.
\end{proof}

\begin{corollary}\label{c:kerQ-Lagrangian}
Let $(V,\omega)$ be a finite-dimensional complex symplectic vector
space with three Lagrangian subspaces $\alpha,\beta,\gamma$. Then we
have $\ker Q(\alpha,\beta;\gamma)=\alpha\cap\beta+\alpha\cap\gamma$.
\end{corollary}
We now study when the form $Q(\alpha,\beta;\gamma)$ is zero if $\ker Q(\alpha,\beta;\gamma)=\alpha\cap\beta+\alpha\cap\gamma$
\begin{lemma}\label{l:3-equivalence}
Let $V$ be an Abelian group with three subgroups
$\alpha,\beta,\gamma$. Then the following three conditions are
equivalent:
\begin{itemize}
  \item[(i)] $ \alpha\cap\beta+\alpha\cap\gamma =(\beta+\gamma)\cap\alpha$,
  \item [(ii)] $\alpha\cap\beta+\beta\cap\alpha = (\alpha+\gamma)\cap\beta$,
  \item [(iii)] $\alpha\cap\gamma+\beta\cap\gamma=(\alpha+\beta)\cap\gamma$.
\end{itemize}
\end{lemma}
\begin{proof}
By the symmetry of the statement, we only need to prove that
(i)$\Rightarrow$(ii).

Assume that (i) holds. Clearly we have
$\alpha\cap\beta+\beta\cap\gamma\subset(\alpha+\gamma)\cap\beta$.
Let $y\in(\alpha+\gamma)\cap\beta$. Then there exist $x\in\alpha$
and $z\in\gamma$ such that $y=x+z$. So
$x=y-z\in(\beta+\gamma)\cap\alpha$. By (i), there exist
$y_1\in\alpha\cap\beta$ and $z_1\in \alpha\cap\gamma$ such that
$x=y_1+z_1$. So $z_1+z=y-y_1\in \beta\cap\gamma$ and
$y=y_1+(z_1+z)\in\alpha\cap\beta+\beta\cap\gamma$.
\end{proof}

\begin{lemma}\label{l:linear-subspaces-dimension}
Let $V$ be a vector space with three finite-dimensional linear
subspaces $\alpha,\beta,\gamma$. Then we have
\begin{equation}\label{e:linear-subspaces-dimension}
\begin{split}
  \dim (\alpha\cap\beta)&+ \dim(\alpha\cap\gamma)+\dim (\beta\cap\gamma) \leq \dim \alpha +\dim \beta + \dim \gamma \\
    &+ \dim (\alpha\cap\beta\cap\gamma)-\dim(\alpha+\beta+\gamma).
\end{split}
\end{equation}
The equality in (\ref{e:linear-subspaces-dimension}) holds if and
only if $\alpha\cap\beta+\alpha\cap\gamma=\alpha\cap(\beta+\gamma)$.
\end{lemma}

\begin{proof}
We have $\dim (\beta\cap\gamma)=\dim \beta +\dim \gamma - \dim
(\beta+\gamma)$. Since $\alpha\cap\beta+\alpha\cap\gamma \subset
\alpha\cap(\beta+\gamma)$, we have
\begin{equation*}
    \begin{split}
      \dim (\alpha\cap\beta)&+\dim (\alpha\cap\gamma)=\dim(\alpha\cap\beta + \alpha\cap\gamma) +\dim(\alpha\cap\beta\cap\gamma)  \\
       &\leq \dim (\alpha\cap(\beta+\gamma)) + \dim(\alpha\cap\beta\cap\gamma)\\
        &=\dim \alpha +\dim (\beta+\gamma)-\dim (\alpha
        +\beta+\gamma)+\dim(\alpha\cap\beta\cap\gamma).
    \end{split}
\end{equation*}
So (\ref{e:linear-subspaces-dimension}) holds. The equality in
(\ref{e:linear-subspaces-dimension}) holds if and only if
$\dim(\alpha\cap\beta+\alpha\cap\gamma)=\dim
(\alpha\cap(\beta+\gamma))$, if and only if $\alpha\cap\beta
+\alpha\cap\gamma=\alpha\cap(\beta+\gamma)$.
\end{proof}

\begin{corollary}\label{c:Lagrangian-subspaces-dimension}
Let $(V,\omega)$ be a complex symplectic vector space of dimension
$2n$ with three Lagrangian subspaces $\alpha,\beta,\gamma$. Then we
have
\begin{equation}\label{e:Lagrangian-subspaces-dimension}
    \dim(\alpha\cap\beta)+\dim(\alpha\cap\gamma)+\dim(\beta\cap\gamma)\leq
    n+2\dim(\alpha\cap\beta\cap\gamma).
\end{equation}
The equality in (\ref{e:Lagrangian-subspaces-dimension}) holds if
and only if $\alpha\cap\beta +
\alpha\cap\gamma=\alpha\cap(\beta+\gamma)$, if and only if
$Q(\alpha,\beta;\gamma)=0$.
\end{corollary}
\begin{proof}
Since $\alpha,\beta,\gamma$ are Lagrangian subspaces of $V$ and
$\dim V=2n$, we have $\dim \alpha=\dim \beta=\dim \gamma=n$ and
\begin{equation*}
    \dim
    (\alpha\cap\beta\cap\gamma)=\dim((\alpha+\beta+\gamma)^{\omega})=2n-\dim(\alpha+\beta+\gamma).
\end{equation*}
By Lemma \ref{l:linear-subspaces-dimension}, the inequality
(\ref{e:Lagrangian-subspaces-dimension}) holds, and the equality in
(\ref{e:Lagrangian-subspaces-dimension}) holds if and only if
$\alpha\cap\beta+\alpha\cap\gamma=\alpha\cap(\beta+\gamma)$. By
Corollary \ref{c:kerQ-Lagrangian},
$\alpha\cap\beta+\alpha\cap\gamma=\alpha\cap(\beta+\gamma)$ holds if
and only if $Q(\alpha,\beta;\gamma)=0$.
\end{proof}

If $\epsilon$ is an isotropic subspace of $(V,\omega)$ such that
$\epsilon=\epsilon^{\omega\omega}$, $\omega$ defines a symplectic
form $\tilde{\omega}$ on $\epsilon^{\omega}/\epsilon$. Moreover, for
each isotropic subspace $\delta$, the image $\pi_{\epsilon}(\delta)$
of $\delta\cap\epsilon^{\omega}$ under the canonical homomorphism:
$\pi_{\epsilon}:\epsilon^{\omega}\rightarrow
\epsilon^{\omega}/\epsilon$ is an isotropic subspace of
$(\epsilon^{\omega}/\epsilon,\tilde{\omega})$.

Let $\alpha,\beta,\gamma$ be three isotropic subspaces of $V$.
Assume that $\epsilon\subset\beta\subset\epsilon^{\omega}$. Then we
have $\beta\cap\epsilon^{\omega}+\epsilon=\beta$, and
\begin{equation*}
    \begin{split}
      (\alpha\cap\epsilon^{\omega}+\epsilon)&\cap(\beta+\gamma\cap\epsilon^{\omega}+\epsilon)=(\alpha+\epsilon)\cap\epsilon^{\omega}\cap(\beta+\gamma)\cap\epsilon^{\omega} \\
         &=(\alpha\cap(\beta+\gamma)+\epsilon)\cap\epsilon^{\omega}\\
        &=\alpha\cap(\beta+\gamma)\cap\epsilon^{\omega}+\epsilon.
    \end{split}
\end{equation*}
So we have
$\pi(\alpha\cap(\beta+\gamma))=(\pi\alpha)\cap(\pi\beta+\pi\gamma)$,
and (cf. \cite[(2.11)]{Du76})
\begin{equation}\label{e:Q-index-pi}
Q(\alpha,\beta;\gamma)(x_1,x_2)=Q(\pi\alpha,\pi\beta;\pi\gamma)(\pi
x_1,\pi x_2),
\end{equation}
here $x_1,x_2\in \alpha\cap(\beta+\gamma)\cap \epsilon^{\omega}$ and
$\pi=\pi_{\epsilon}$.

\begin{lemma}\label{l:Q-definite-path}
Let $(V,\omega)$ be a finite-dimensional complex symplectic vector
space. Let $\alpha(s)$, $s\in(-\varepsilon,\varepsilon)$ be a path
in $\mathcal{L}(V)$ differentiable at $s=0$. Let $\beta$ be a
Lagrangian subspace of $V$. Assume that $Q(\alpha,0)$ is positive
definite. Then there exists an $\varepsilon_1\in(0,\epsilon)$ such that for
$s\in(0,\varepsilon_1)$, we have $V=\alpha(\pm s)\oplus \beta$, and
\begin{equation}\label{e:Q-definite-path}
    m^{\mp}(Q(\alpha(0),\beta;\alpha(\pm s)))=m^{\mp}(Q(\beta,\alpha(\pm
    s);\alpha(0)))=0.
\end{equation}
\end{lemma}
\begin{proof}
By Proposition \ref{p:calculate-maslov-index}, there exists an
$\varepsilon_2\in(0,\epsilon)$ such that $V=\alpha(\pm s)\oplus
\beta$ holds for each $s\in (0,\varepsilon_2)$. Set
$\epsilon:=\alpha(0)\cap\beta$ and $\pi:=\pi_{\epsilon}$. By the
proof of \cite[Corollary 1.3.4]{BoZh14}, there is a
$\tilde{\beta}\in \mathcal{L}(V)$ such that
$V=\alpha(0)\oplus\tilde{\beta}$ and $\pi \tilde{\beta}=\pi \beta$.
Since $Q(\alpha,0)$ is positive definite, there exists an
$\varepsilon_1\in(0,\varepsilon_2)$ such that for
$s\in(0,\varepsilon_1)$, the form
$Q(\alpha(0),\tilde{\beta};\alpha(\pm s))$ is positive (negative)
definite. Let $s$ be in $(0,\varepsilon_1)$. By
(\ref{e:Q-index-pi}), the form $Q(\pi \alpha(0),\pi
\tilde{\beta};\pi \alpha(\pm s))$ is positively (negatively)
definite. Note that $\alpha(0)\subset \epsilon^{\omega}$. By
(\ref{e:Q-index-pi}) again, the form $Q(\alpha(0),\beta;\alpha(\pm
s))$ is positive (negative) definite. By Lemma \ref{l:Q-cyclic},
(\ref{e:Q-definite-path}) holds.
\end{proof}

\subsection{The triple index and the H{\"o}rmander index}\label{ss:triple-and-Hormander-index}
Let $(V,\omega)$ be a complex symplectic vector space of dimension
$2n$. Let $\lambda_1,\lambda_2,\mu_1,\mu_2$ be four Lagrangian
subspaces of $(V,\omega)$.

\begin{definition}\label{d:Hormander-index}(\cite[Definition 3.4.4]{BoZh14})
Assume that there are continuous paths $\lambda(s)$ and $\mu(s)$,
$s\in[a,b]$ of Lagrangian subspaces of $(V,\omega)$ such that
$\lambda(a)=\lambda_1$, $\lambda(b)=\lambda_2$, $\mu(a)=\mu_1$,
$\mu(b)=\mu_2$. Then the \textit{H{\"o}rmander index}
$s(\lambda_1,\lambda_2;\mu_1,\mu_2)$ is defined by
\begin{eqnarray}
% \nonumber to remove numbering (before each equation)
  s(\lambda_1,\lambda_2;\mu_1,\mu_2) &=& \Mas \{\lambda,\mu_2\}-\Mas\{\lambda,\mu_1\}\label{e:path-front} \\
   &=&\Mas\{\lambda_2,\mu\}-\Mas\{\lambda_1,\mu\}.\label{e:path-back}
\end{eqnarray}
\end{definition}

By \cite[Proposition 2.3.1.b,f]{BoZh14}, we have
\begin{eqnarray}
% \nonumber to remove numbering (before each equation)
  &s(\lambda_1,&\lambda_2;\mu_1,\mu_3) = s(\lambda_1,\lambda_2;\mu_1,\mu_2)+s(\lambda_1,\lambda_2;\mu_2,\mu_3),\label{e:s-adjacent}\\
  &s(\lambda_1,&\lambda_2;\mu_1,\mu_2) = -s(\lambda_1,\lambda_2;\mu_2,\mu_1)\label{e:s-opposite}\\
   &&=-s(\mu_1,\mu_2;\lambda_1,\lambda_2)+\sum_{j,k\in\{1,2\}}(-1)^{j+k+1}
   \dim (\lambda_j\cap\mu_k).\label{e:s-front-back}
\end{eqnarray}

\begin{lemma}\label{l:partial-transversal-case}
Let $\lambda(s)$, $s\in[a,b]$ be a Lagrangian path of complex
symplectic vector space $(V,\omega)$. Let $\alpha\in\mathcal{L}(V)$
be such that $\alpha\cap\lambda(a)=\alpha\cap\lambda(b)=\{0\}$. Then
we have
\begin{eqnarray}
% \nonumber to remove numbering (before each equation)
  s(\alpha,\lambda(a);\lambda(a),\lambda(b)) &=& m^-(Q(\lambda(a),\alpha;\lambda(b))),\label{e:partial-transversal-case-I} \\
  \Mas\{\lambda(a),\lambda\} &=&
  \Mas\{\alpha,\lambda\}+m^-(Q(\lambda(a),\alpha;\lambda(b))).\label{e:partial-transversal-case-II}
\end{eqnarray}
\end{lemma}
\begin{proof}
Let $A(s)\in \Hom (\lambda(a),\alpha)$, $s\in[a,b]$ be a path of
linear maps such that, the form $\omega(x,A(s)y)$,
$x,y\in\lambda(a)$ is Hermitian for each $s\in[a,b]$, $A(a)=0$,
$\lambda(b)=\Graph(A(b))$. Consider a special path
$\lambda(s)=\Graph(A(s))$. By \cite[Lemma 2.3.2]{BoZh14}, we have
\begin{equation*}
    \begin{split}
      s(\alpha,\lambda(a);\lambda(a),\lambda(b)) & =\Mas\{\lambda(a),\lambda\}-\Mas\{\alpha,\lambda\} \\
        &=-\Mas_{-}\{\lambda,\lambda(a)\}+\Mas_{-}\{\lambda,\alpha\}\\
        &=m^-(Q(\lambda(a),\alpha;\lambda(b))).
    \end{split}
\end{equation*}
The equality (\ref{e:partial-transversal-case-II}) follows from
(\ref{e:partial-transversal-case-I}).
\end{proof}

\begin{corollary}\label{c:transversal-case}(\cite[(3.3.5),(3.3.7)]{Ho71}),
\cite[(2.10),(2.13)]{Du76}) Let $\lambda_1,\lambda_2,\mu_1,\mu_2$ be
Lagrangian subspaces of complex symplectic vector space
$(V,\omega)$. Assume that $\lambda_i\cap\mu_j=0$, $i,j=1,2$. Then we
have
\begin{eqnarray}
% \nonumber to remove numbering (before each equation)
  s(\lambda_1,\lambda_2;\mu_1,\mu_2) &=& -s(\mu_1,\mu_2;\lambda_1,\lambda_2)\label{e:transversal-case-I} \\
   &=&
   m^-(Q(\lambda_1,\lambda_2;\mu_1))-m^-(Q(\lambda_1,\lambda_2;\mu_2)).
\end{eqnarray}
\end{corollary}
\begin{proof}
By (\ref{e:s-front-back}), (\ref{e:s-adjacent}),
(\ref{e:s-opposite}), Lemma \ref{l:partial-transversal-case},
Corollary \ref{c:kerQ-Lagrangian} and (\ref{e:Q-transposition}) we
have
\begin{equation*}
    \begin{split}
       s(\lambda_1,\lambda_2;\mu_1,\mu_2) & =-s(\mu_1,\mu_2;\lambda_1,\lambda_2) \\
         &=s(\mu_2,\lambda_1;\lambda_1,\lambda_2)-s(\mu_1,\lambda_1;\lambda_1,\lambda_2)\\
         &=m^-(Q(\lambda_1,\mu_2;\lambda_2))-m^-(Q(\lambda_1,\mu_1;\lambda_2))\\
         &=m^-(Q(\lambda_1,\lambda_2;\mu_1))-m^-(Q(\lambda_1,\lambda_2;\mu_2)).
\end{split}
\end{equation*}
\end{proof}

The following corollary was proved by J. J. Duistermaat
\cite[(2.16)]{Du76}.

\begin{corollary}\label{c:i-index-definition}
Let $\alpha,\beta,\gamma$ be Lagrangian subspaces of complex
symplectic vector space $V$. Define the triple index of
$\alpha,\beta,\gamma$ by
\begin{equation}\label{e:i-index-definition}
i(\alpha,\beta,\gamma)=m^-(Q(\alpha,\delta;\beta))+m^-(Q(\beta,\delta;\gamma))-m^-(Q(\alpha,\delta;\gamma)),
\end{equation}
where $\delta\in\mathcal{L}(V)$ be such that
$\delta\cap\alpha=\delta\cap\beta=\delta\cap\gamma=\{0\}$. Then the
triple index is well-defined.
\end{corollary}

Now we calculate the triple index $i(\alpha,\beta,\gamma)$.

\begin{lemma}\label{l:calculate-triple-index}
Let $\alpha,\beta,\gamma$ be three Lagrangian subspaces of $V$. Then
we have
\begin{eqnarray}
% \nonumber to remove numbering (before each equation)
  i(\alpha,\beta,\gamma) &=& m^+(Q(\alpha,\beta;\gamma))+\dim (\alpha\cap\gamma)-\dim (\alpha\cap\beta\cap\gamma) \label{e:equality-triple-index}\\
   &\leq& n-\dim (\alpha\cap\beta)-\dim(\beta\cap\gamma)+\dim
   (\alpha\cap\beta\cap\gamma).\label{e:inequality-triple-index}
\end{eqnarray}
\end{lemma}
\begin{proof}
Denote by $\epsilon:=\alpha\cap\beta+\beta\cap\gamma$ and
$\pi:=\pi_{\epsilon}$. Recall that
$\pi\alpha=(\alpha+\epsilon)\cap\epsilon^{\omega}/\epsilon$. Note
that $\epsilon\subset\beta$,
\begin{eqnarray*}
% \nonumber to remove numbering (before each equation)
  \epsilon^{\omega} &=& (\alpha+\beta)\cap(\beta+\gamma),\ and \\
  \alpha+\epsilon &=&
  \alpha+\alpha\cap\beta+\beta\cap\gamma=\alpha+\beta\cap\gamma.
\end{eqnarray*}
It follows that $(\pi\alpha)\cap(\pi\beta)=\{0\}$, and
\begin{equation*}\begin{split}
                  (\pi\alpha)\cap(\pi\gamma) & =((\alpha+\epsilon)\cap(\alpha+\beta)\cap(\beta+\gamma)\cap(\gamma+\epsilon))/\epsilon \\
                    &
                    =((\alpha+\epsilon)\cap(\gamma+\epsilon))/\epsilon\\
                    &=((\alpha+\epsilon)\cap\gamma+\epsilon)/\epsilon\\
                    &=((\alpha+\beta\cap\gamma)\cap\gamma+\epsilon)/\epsilon\\
                    &=(\alpha\cap\gamma+\beta\cap\gamma+\alpha\cap\beta+\beta\cap\gamma)/(\alpha\cap\beta+\beta\cap\gamma)\\
                    &\cong(\alpha\cap\gamma)/(\alpha\cap\gamma\cap(\alpha\cap\beta+\beta\cap\gamma))\\
                    &=(\alpha\cap\gamma)/(\alpha\cap\beta\cap\gamma).
                \end{split}
\end{equation*}
Then by Lemma \ref{l:kerQ-isotopic}, we have
\begin{equation*}
    \begin{split}
      m^0(Q(\pi\alpha,\pi\beta;\pi\gamma)) & = \dim((\pi\alpha)\cap(\pi\beta))+\dim((\pi\alpha)\cap(\pi\gamma))\\
        & =\dim (\alpha\cap\gamma)-\dim(\alpha\cap\beta\cap\gamma).
    \end{split}
\end{equation*}

By Lemma \ref{l:Q-cyclic}, (\ref{e:Q-transposition}), Remark
\ref{r:Q-Duistermaat} and \cite[Lemma 2.4]{Du76}, we have
\begin{equation*}
    \begin{split}
      i(\alpha,\beta,\gamma) & =(m^++m^0)(Q(\pi \alpha,\pi \beta;\pi \gamma)) \\
        & =m^+(Q(\alpha,\beta;\gamma))+\dim
        (\alpha\cap\gamma)-\dim(\alpha\cap\beta\cap\gamma)\\
        &\leq\dim(\pi\alpha)=n-\dim \epsilon\\
        &=n-\dim(\alpha\cap\beta)-\dim(\beta\cap\gamma)+\dim(\alpha\cap\beta\cap\gamma).
    \end{split}
\end{equation*}
\end{proof}

Let $\alpha,\beta$ be two Lagrangian subspaces of $V$. We define
\begin{equation}\label{e:definition-L0}
\mathcal{L}_0(\alpha,\beta):=\{\gamma\in\mathcal{L}(V);\alpha\cap\gamma+\beta\cap\gamma=(\alpha+\beta)\cap\gamma\}.
\end{equation}

\begin{corollary}\label{c:kerQ-0}
Let $\alpha,\beta,\gamma$ be three Lagrangian subspaces of $V$.
Assume that $\gamma\in\mathcal{L}_0(\alpha,\beta)$. Let
$(\alpha_1,\beta_1,\gamma_1)$ be a permutation of
$(\alpha,\beta,\gamma)$. Then we have
\begin{equation}\label{e:i-permutation}
i(\alpha_1,\beta_1,\gamma_1)=\dim (\alpha_1\cap\gamma_1)-\dim
(\alpha_1\cap\beta_1\cap\gamma_1).
\end{equation}
In particular, we have
\begin{equation}\label{e:i-two-repeat}
i(\alpha,\alpha,\beta)=i(\beta,\alpha,\alpha)=0,\quad
i(\alpha,\beta,\alpha)=n-\dim(\alpha\cap\beta).
\end{equation}
\end{corollary}
\begin{proof}
By Corollary \ref{c:kerQ-Lagrangian} and Lemma
\ref{l:3-equivalence}, we have $Q(\alpha_1,\beta_1;\gamma_1)=0$. By
Lemma \ref{l:calculate-triple-index}, (\ref{e:i-permutation}) holds.
Since $\beta\in\mathcal{L}_0(\alpha,\alpha)$, (\ref{e:i-two-repeat})
holds.
\end{proof}

Denote by $J\in \End (V)$ with $\omega(x,y)=\langle Jx,y\rangle$ for
each $x,y\in V$.

By Lemma \ref{l:Q-definite-path}, there is an $\epsilon >0$ such
that, for each $s\in(0,\epsilon)$, we have
$(e^{Js}\gamma)\cap\alpha=(e^{Js}\gamma)\cap\beta=(e^{Js}\gamma)\cap\gamma=\{0\}$,
and
\begin{equation*}
    m^-(Q(\beta,e^{Js}\gamma;\gamma))=m^-(Q(\alpha,e^{Js}\gamma;\gamma))=0.
\end{equation*}
By Corollary \ref{c:i-index-definition}, we have
\begin{equation}\label{e:i-positive-path}
i(\alpha,\beta,\gamma)=m^-(Q(\alpha,e^{Js}\gamma;\beta)).
\end{equation}

\begin{lemma}\label{l:s-special-path}
Let $\lambda,\mu:[a,b]\rightarrow\mathcal{L}(V)$ be two paths of
Lagrangian subspaces of $V$. Then there is an $\epsilon>0$ such that
for each $s\in(0,\epsilon)$, we have
\begin{equation}\label{e:s-special-path}
s(\lambda(a),\lambda(b);e^{Js}\mu(a),e^{Js}\mu(b))=s(\lambda(a),\lambda(b);\mu(a),\mu(b)).
\end{equation}
\end{lemma}
\begin{proof}
Let $s>0$ be sufficiently small. According to Proposition
\ref{p:calculate-maslov-index}, we have $\Mas
\{\lambda(s_1),e^{Jt}\mu(s_2);t\in[0,s]\}=0$, for
$s_1,s_2\in\{a,b\}$. So by the definition of H{\"o}rmander index
(\ref{e:path-back}), we have
$s(\lambda(a),\lambda(b);\mu(s_1),e^{Js}\mu(s_1))=0$, for $s_1=a,b$.
Meanwhile by (\ref{e:s-adjacent}) and (\ref{e:s-opposite}),
(\ref{e:s-special-path}) holds.
\end{proof}

\begin{proof}[Proof of Theorem \ref{t:calculate-hormander-index}]
Let $s>0$ be sufficiently small. By Lemma \ref{l:s-special-path},
the proof of Corollary \ref{c:transversal-case} and
(\ref{e:i-positive-path}) we have
\begin{equation*}
\begin{split}
  s(\lambda_1,\lambda_2&;\mu_1,\mu_2)= s(\lambda_1,\lambda_2;e^{Js}\mu_1,e^{Js}\mu_2)\\
    &=m^-(Q(\lambda_1,e^{Js}\mu_2;\lambda_2))-m^-(Q(\lambda_1,e^{Js}\mu_1;\lambda_2))\\
    &=i(\lambda_1,\lambda_2,\mu_2)-i(\lambda_1,\lambda_2,\mu_1).
\end{split}
\end{equation*}

By (\ref{e:s-front-back}), Lemma \ref{l:calculate-triple-index} and
Lemma \ref{l:Q-cyclic}, we have
\begin{equation*}
    \begin{split}
       s(\lambda_1,\lambda_2&;\mu_1,\mu_2)=-s(\mu_1,\mu_2;\lambda_1,\lambda_2)+\sum_{j,k\in\{1,2\}}(-1)^{j+k+1}\dim(\lambda_j\cap\mu_k)  \\
         &=i(\mu_1,\mu_2,\lambda_1)-i(\mu_1,\mu_2,\lambda_2)+\sum_{j,k\in\{1,2\}}(-1)^{j+k+1}\dim(\lambda_j\cap\mu_k)\\
         &=m^+(Q(\mu_1,\mu_2;\lambda_1))+\dim(\lambda_1\cap\mu_1)-\dim(\lambda_1\cap\mu_1\cap\mu_2)\\
         &\quad
         -m^+(Q(\mu_1,\mu_2;\lambda_2))-\dim(\lambda_2\cap\mu_1)+\dim(\lambda_2\cap\mu_1\cap\mu_2)\\
         &\quad
         +\sum_{j,k\in\{1,2\}}(-1)^{j+k+1}\dim(\lambda_j\cap\mu_k)\\
         &=m^+(Q(\lambda_1,\mu_1;\mu_2))+\dim(\lambda_1\cap\mu_2)-\dim(\lambda_1\cap\mu_1\cap\mu_2)\\
         &\quad
         -m^+(Q(\lambda_2,\mu_1;\mu_2))-\dim(\lambda_2\cap\mu_2)+\dim(\lambda_2\cap\mu_1\cap\mu_2)\\
         &=i(\lambda_1,\mu_1,\mu_2)-i(\lambda_2,\mu_1,\mu_2).
         \end{split}
\end{equation*}
\end{proof}

\begin{corollary}\label{c:s-2-repeat}
Let $(V,\omega)$ be a complex symplectic vector space of dimension
$2n$. Let $\lambda\in C([a,b],\mathcal{L}(V))$ be a Lagrangian path.
Then for each $\mu\in\mathcal{L}(V)$, we have
\begin{eqnarray}
    &s(\lambda(a),&\lambda(b);\lambda(a),\mu)=-i(\lambda(b),\lambda(a),\mu)\label{e:s-13-1}\\
    &&\leq\dim(\lambda(a)\cap\lambda(b)\cap\mu)-\dim(\lambda(b)\cap\mu)\leq0,\label{e:s-13-2}\\
&s(\lambda(a),&\lambda(b);\lambda(a),\mu)\geq
\dim(\lambda(a)\cap\lambda(b))+\dim(\lambda(a)\cap\mu) \label{e:s-13-3}\\
&&\mspace{125mu} -\dim(\lambda(a)\cap\lambda(b)\cap\mu)-n,\nonumber\\
&s(\lambda(a),&\lambda(b);\lambda(b),\mu)=i(\lambda(a),\lambda(b),\mu)\geq0,\label{e:s-23}\\
 &&\Mas\{\lambda,\lambda(b)\}\leq\Mas\{\lambda,\mu\}\leq\Mas\{\lambda,\lambda(a)\}.\label{e:Maslov-index-inequality}
\end{eqnarray}
\end{corollary}
\begin{proof}
By Theorem \ref{t:calculate-hormander-index}, Lemma
\ref{l:calculate-triple-index} and (\ref{e:i-two-repeat}), we have
\begin{equation*}
    \begin{split}
      s(\lambda(a),\lambda(b);\lambda(a),\mu) & =i(\lambda(a),\lambda(a),\mu)-i(\lambda(b),\lambda(a),\mu) \\
        & =-i(\lambda(b),\lambda(a),\mu)\\
        &\leq\dim
        (\lambda(a)\cap\lambda(b)\cap\mu)-\dim(\lambda(b)\cap\mu)\leq0,\\
        s(\lambda(a),\lambda(b);\lambda(a),\mu)&\geq
        \dim(\lambda(a)\cap\lambda(b))+\dim (\lambda(a)\cap\mu)\\
        &\quad -\dim(\lambda(a)\cap\lambda(b)\cap\mu)-n,\\
        s(\lambda(a),\lambda(b);\lambda(b),\mu)&=i(\lambda(a),\lambda(b),\mu)-i(\lambda(a),\lambda(b),\lambda(b))\\
        &=i(\lambda(a),\lambda(b),\mu)\geq0.
    \end{split}
\end{equation*}
By Definition \ref{d:Hormander-index} we obtain
(\ref{e:Maslov-index-inequality}).
\end{proof}

\section{Iteration inequalities of Maslov-type index}\label{s:iteration-inequalities}
\subsection{Iteration inequalities with periodical boundary condition}\label{ss:iteration-inequalities-periodical}
Let $\gamma_l\in\mathcal{P}_{\tau_l}(V)$, $l=1,\ldots,k$ be $k$
symplectic paths starting from the identity $I_{V}$, where
$\tau_l\geq0$. For each $l=1,\ldots,k$, we set
\begin{equation}\label{e:definition-matrix}
M_l:=\gamma_l(\tau_l),\quad T_l:=\sum_{j=1}^l\tau_j,\quad
\tilde{M_l}=\prod_{j=1}^l M_{l+1-j}.
\end{equation}
The \textit{iteration} of $\gamma_1,...,\gamma_k$ is a symplectic
path $\tilde{\gamma}\in\mathcal{P}_{T_k}(V)$ defined by
\begin{equation}\label{e:definition-iteration-path}
    \tilde{\gamma}(t)=\gamma_l(t-T_{l-1})\tilde{M}_{l-1},\quad
    T_{l-1}\leq t\leq T_l,\ l=1,...k.
\end{equation}

Let $\tau>0$ be a positive number and $k\in\NN$ be a positive
integer. For each symplectic path $\gamma\in\mathcal{P}_{\tau}(V)$,
we define its \textit{$k$-th} \textit{iteration} $(\gamma,k)$ to be
the iteration path of $k$ copies of $\gamma$.

If $\gamma$ is a real symplectic path, the iteration inequalities
for $(\gamma,k)$ with periodical boundary condition was obtained by
C. Liu and Y. Long in 1997 (cf. \cite[Theorem 10.1.3]{Lo02}).

Our iteration inequalities with periodic boundary condition read as
follows.

\begin{theorem}\label{t:itertion-inequalities-any}
Let $(V,\omega)$ be a symplectic vector space of dimension $m$ and
$k\geq2$ be an integer. Let $\gamma_l\in\mathcal{P}_{\tau_l}(V)$,
$l=1,...,k$ be $k$ symplectic paths starting from the identity
$I_V$, where $\tau_l\geq0$. Let
$\tilde{\gamma}\in\mathcal{P}_{T_k}(\gamma)$ be the iteration of
$\gamma_1,...,\gamma_k$. Let $M_l$ and $\tilde{M}_l$ be defined by
$(\ref{e:definition-matrix})$. Denote by $\nu_1(M):=\dim \ker
(M-I_V)$ and $\mathfrak{N}_1(M):=\ker(M-I_V)$ for $M\in\Sp(V)$. For
$\mathcal{M}:=(M_1,M_2,...,M_k)$, we define
\begin{eqnarray}
% \nonumber to remove numbering (before each equation)
  A(\mathcal{M}) &:=& \dim(\bigcap_{l=1}^k\mathfrak{N}_1(M_l)),\label{e:A-M} \\
  B(\mathcal{M})&:=& \sum_{l=2}^k \dim(\mathfrak{N}_1(M_l)\cap\mathfrak{N}_1(\tilde{M}_{l-1}))-\sum_{l=2}^{k-1}\nu_1(\tilde{M}_l).
\label{e:B-M}
\end{eqnarray}
Then we have
\begin{eqnarray}
% \nonumber to remove numbering (before each equation)
  \sum_{l=1}^k\nu_1(M_l)-B(\mathcal{M})-m(k-1) &\leq& i_1(\tilde{\gamma})-\sum_{l=1}^ki_1(\gamma_l)\label{e:-B-Maslov-type} \\
  &\leq& B(\mathcal{M})-\nu_1(\tilde{M}_k),\label{e:B-} \\
  B(\mathcal{M}) &\leq& A(\mathcal{M}).\label{e:B-A}
\end{eqnarray}
 \end{theorem}
\begin{proof}
We divide the proof into four steps.

\textbf{Step 1.} The case that $k=2$.

In this case we have $A(\mathcal{M})=B(\mathcal{M})$. By Lemma
\ref{l:calculate-Maslov-multiply}, the definition of the Maslov-type
index and the H{\"o}rmander index we have
\begin{equation*}
    \begin{split}
       i_1(\tilde{\gamma})&-i_1(\gamma_1)-i_1(\gamma_2)=i_1(\gamma_2M_1)-i_1(\gamma_2)  \\
         &=\Mas\{\gamma_2,\Graph(M_1^{-1})\}-\Mas\{\gamma_2,\Graph(I_V)\}\\
         &=s(\Graph(I_V),\Graph(M_2);\Graph(I_V),\Graph(M_1^{-1})).
     \end{split}
\end{equation*}
By Corollary \ref{c:s-2-repeat} we have
\begin{equation}\label{e:iteration-inequality-two}
    \begin{split}
      \nu_1(M_1)&+\nu_1(M_2)- \dim(\ker(M_2-I_V)\cap\ker(M_1-I_V))-m \\
        & \leq i_1(\tilde{\gamma})-i_1(\gamma_1)-i_1(\gamma_2)\\
        &= -i(\Graph(M_2),\Graph(I_V),\Graph(M_1^{-1}))\\
        &\leq -\nu_1(M_1M_2)+\dim (\ker(M_2-I_V)\cap\ker(M_1-I_V)).
    \end{split}
\end{equation}
So the case that $k=2$ follows.

\textbf{Step 2.} The inequalities (\ref{e:-B-Maslov-type}) and
(\ref{e:B-}) hold.

By Step 1, for $l=2,...,k$ we have
\begin{equation}\label{e:iteration-inequality-l}
    \begin{split}
      \nu_1(\tilde{M}_{l-1})&+\nu_1(M_l)-\dim \mathfrak{N}_1(M_l)\cap\mathfrak{N}_1(\tilde{M}_{l-1})-m\\
        & \leq
        i_1(\tilde{\gamma}|_{[0,T_l]})-i_1(\tilde{\gamma}|_{[0,T_{l-1}]})-i_1(\gamma_l)\\
        & \leq -\nu_1(\tilde{M}_l) +\dim
        \mathfrak{N}_1(M_l)\cap\mathfrak{N}_1(\tilde{M}_{l-1}).
    \end{split}
\end{equation}
Add up (\ref{e:iteration-inequality-l}) for $l=2,...,k$, we obtain
(\ref{e:-B-Maslov-type}) and (\ref{e:B-}).

\textbf{Step 3.} Let $V_1,V_2,V_3$ be linear subspaces of a vector
space $V$. Assume that $V_1\subset V_2$ and $\dim V_2<+\infty$. Then
we have
\begin{equation*}
    \frac{V_1}{V_1\cap V_3}\cong \frac{V_1+V_3}{V_3} \subset
    \frac{V_2+V_3}{V_3}\cong \frac{V_2}{V_2\cap V_3}.
\end{equation*}
So we obtain
\begin{equation}\label{e:dim-V_1}
    \dim \frac{V_1}{V_1\cap V_3}\leq \dim \frac{V_2}{V_2\cap V_3}.
\end{equation}

\textbf{Step 4.} The inequality (\ref{e:B-A}) holds.

For $l=2,...,k$ we have $\cap_{j=1}^{l-1}\mathfrak{N}(M_j)\subset
\mathfrak{N}(\tilde{M}_{l-1})$. By Step 3 we have
\begin{equation}\label{e:dim-ker-cap}
\dim
\frac{\cap_{j=1}^{l-1}\mathfrak{N}(M_j)}{\cap_{j=1}^l\mathfrak{N}(M_j)}\leq
\dim
\frac{\mathfrak{N}(\tilde{M}_{l-1})}{\mathfrak{N}(\tilde{M}_{l-1})\cap
\mathfrak{N}(M_l)}.
\end{equation}
Add up (\ref{e:dim-ker-cap}) for $l=2,...,k$, we obtain
\begin{equation*}
    \nu_1(M_1)-A(\mathcal{M})\leq \nu_1(M_1)-B(\mathcal{M}).
\end{equation*}
So (\ref{e:B-A}) holds.
\end{proof}

Now we generalize \cite[Thereom 10.2.2]{Lo02} to the complex case.
Our method gives another proof of \cite[Theorem 10.2.2]{Lo02}.

\begin{rem}\label{r:long-ineq}
%The normal forms of the real symplectic matrix in Y. Long's book \cite{Lo02} do the first author a great favor.
Note that Y. Long's book \cite{Lo02} deals with standard real
symplectic space $(\RR^{2n},\omega)$, the real symplectic matrix
$M\in\Sp(2n,\RR):=\{M\in\RR^{2n\times2n};M^{T}JM=J\}$, where $M^T$ is the transpose of matrix $M$ and
$J=\left(
                                                                    \begin{array}{cc}
                                                                      0 & -I_n\\
                                                                      I_n & 0 \\
                                                                    \end{array}
                                                                  \right)
$. The set of real symplectic paths from the identity is defined by
\cite[2.0.1]{Lo02}:
$$\mathcal{P}_{\tau}(2n,\RR):=\{\gamma\in
C([0,\tau],\Sp(2n,\RR));\gamma(0)=I_{2n}\}.$$ Let
$(\mathbb{C}^{2n},\omega)$ be the standard complex symplectic space,
where $\omega(x,y)=\langle Jx,y\rangle$, $\forall x,y \in
\mathbb{C}^{2n}$. Since
$\CC^{2n}=\RR^{2n}\otimes\CC=\RR^{2n}\oplus\sqrt{-1}\RR^{2n}$,
$\Sp(2n,\RR)\subset \Sp(\CC^{2n})$.Thus if $\gamma\in
\mathcal{P}_{\tau}(2n,\RR)$, then
$\gamma\in\mathcal{P}_{\tau}(\CC^{2n})$, so $i_1(\gamma)$ has been
defined in Definition \ref{d:Maslov-type-index}.

We have these facts : if a path $\gamma\in\mathcal{P}_{\tau}(2n,\RR)$,
then
\begin{eqnarray*}
% \nonumber to remove numbering (before each equation)
  i_1(\gamma)\   &=& i_1^L(\gamma)\ (\text{in \cite[Chapter 5]{Lo02}}) +n,\\
  \nu_1(M)\   &=& \nu_1^L(M)\ (\text{in \cite[Definition 5.1]{Lo02}}).
\end{eqnarray*}
\end{rem}

For $M\in\Sp(V)$, we define the \textit{elliptic height} $e(M)$ of
$M$ by the total algebraic multiplicity of all eigenvalues of $M$ on
the unit circle $S^1$.

\begin{theorem}\label{t:iteration-inequalities}(cf. \cite[Theorem
10.2.2]{Lo02}) Let $(V,\omega)$ be a symplectic vector space of
dimension $m$. For any $\gamma \in \mathcal{P}_{\tau}(V)$, set
$M:=\gamma(\tau)$. Then for any $k_1,k_2\in\NN$, we have
\begin{equation}\label{e:iteration-inequalities}
\begin{split}
  \nu_1(M^{k_1})&+ \nu_1(M^{k_2})-\nu_1(M^{(k_1,k_2)})-\frac{e(M)+m}{2} \\
    & \leq i_1(\gamma,k_1+k_2)-i_1(\gamma,k_1)-i_1(\gamma,k_2)\\
    & \leq \nu_1(M^{(k_1,k_2)})-\nu_1(M^{k_1+k_2})+\frac{e(M)-m}{2},
\end{split}
\end{equation}
where $(k_1,k_2)$ is the greatest common divisor of $k_1$ and $k_2$.
\end{theorem}

\begin{proof}
For each $\lambda\in\CC$, we denote by $E_{\lambda}(M)$ the root
vector space of $M$ belonging to the eigenvalue $\lambda$. Set
\begin{eqnarray*}
% \nonumber to remove numbering (before each equation)
  V_1:=\oplus_{|\lambda|=1}E_{\lambda}(M),&&\quad V_2:=\oplus_{|\lambda|\neq1}E_{\lambda}(M),\\
M_1:=M|_{V_1},&&\quad  M_2:=M|_{V_2},\\
    \alpha:=\oplus_{|\lambda|>1}E_{\lambda}(M),&&\quad \beta:=\oplus_{|\lambda|<1}E_{\lambda}(M),\\
    A:=M|_{\alpha},&&\quad  B:=M|_{\beta}.
\end{eqnarray*}
By a similar proof in \cite[Therorem 1.3.1, Lemma 1.3.2, Theorem
1.7.3]{Lo02}, $V_1,V_2$ are symplectic subspaces of $V$,
$V_1=V_2^{\omega}$, $M_l\in\Sp(V_l)$, $l=1,2$, $V=V_1\oplus V_2$,
$M=M_1\oplus M_2$, $e(M_1)=\dim V_1$, $e(M_2)=0$, $\alpha,\beta $
are Lagrangian subspaces of $V_2$, $V_2=\alpha\oplus\beta$, and
$M_2=A\oplus B$. We denote by $m_l:=\dim V_l$, $l=1,2$. Then we have
$m_2=2\dim\alpha$. Clearly for each $k\in\NN$ we have
\begin{alignat*}{2}
    \nu_1(M^k)&=\nu_1(M_1^k),&\qquad \nu_1(M_2^k)&=0,\\
    e(M)&=e(M_1),&\qquad e(M_2)&=0.
\end{alignat*}

Note that $\ker(M_1^{k_1}-I_{V_1})\cap\ker(M_1^{k_2}-I_{V_1})=\ker
(M_1^{(k_1,k_2)}-I_{V_1})$. By (\ref{e:iteration-inequality-two}) we
have
\begin{equation}\label{e:elliptic case}
\begin{split}
\nu_1(M_1^{k_1})&+ \nu_1(M_1^{k_2})-\nu_1(M_1^{(k_1,k_2)})-m_1 \\
    & \leq -i(\Graph(M_1^{k_2}),\Graph(I_{V_1}),\Graph(M_1^{-k_1}))\\
    & \leq \nu_1(M_1^{(k_1,k_2)})-\nu_1(M_1^{k_1+k_2}).
\end{split}
\end{equation}

Since $e(M_2)=0$, we have
\begin{equation*}
    \Graph(M_2^{k_2})\cap\Graph(I_{V_2})=\Graph(M_2^{-k_1})\cap\Graph(I_{V_2})=\Graph(M_2^{k_2})\cap\Graph(M_2^{-k_1})=\{0\}.
\end{equation*}
So the form
$Q:=Q(\Graph(M_2^{k_2}),\Graph(I_{V_2});\Graph(M_2^{-k_1}))$ is a
non-degenerate form on $\Graph(M_2^{k_2})$. Since
$\alpha\times\alpha\in\mathcal{L}(V_2\times V_2)$ and $M_2$ is in
quasi-diagonal form, that is, with respect a suitable basis for $V_2=\alpha\oplus\beta$, $M_2=\left(\begin{array}{cc}
A & 0\\
0 & (\bar{A}^T)^{-1}\\
\end{array}
\right)
$, $A\in \GL(\CC^{\frac{m_2}{2}})$, we have
\begin{align*}
   & Q|_{\Graph(M_2^{k_2})\cap(\alpha\times\alpha)}=Q|_{\Graph(M_2^{k_2})\cap(\beta\times\beta)}=0,\\
    &\Graph(M_2^{k_2})=(\Graph(M_2^{k_2})\cap(\alpha\times\alpha))\oplus(\Graph(M_2^{k_2})\cap(\beta\times\beta)).
\end{align*}
So we have
\begin{equation}\label{e:hyperbolic-case}
    i(\Graph(M_2^{k_2}),\Graph(I_{V_2}),\Graph(M_2^{-k_1}))=m^+(Q)=\dim\alpha=\frac{m_2}{2}.
\end{equation}

Clearly we have
\begin{equation*}
    \begin{split}
         i(\Graph(M^{k_2}),&\Graph(I_{V}),\Graph(M^{-k_1}))= i(\Graph(M_1^{k_2}),\Graph(I_{V_1}),\Graph(M_1^{-k_1})) \\
         & + i(\Graph(M_2^{k_2}),\Graph(I_{V_2}),\Graph(M_2^{-k_1})).
     \end{split}
\end{equation*}
Our result then follows from (\ref{e:iteration-inequality-two}),
(\ref{e:elliptic case}) and (\ref{e:hyperbolic-case}).
\end{proof}

\subsection{Iteration inequalities with non-periodical boundary
conditions}\label{ss:iteration-inequalities-non-periodical} In this
subsection we derive some iteration inequalities with non-periodical
boundary conditions from Theorem \ref{t:calculate-hormander-index}.

Firstly we give the following useful lemma in the study of
Maslov-type index with brake symmetry.

Let $(V,\omega)$ be a complex symplectic vector space of dimension
$m$. Denote by $(\tilde{V},\tilde{\omega}):=(V\times
V,(-\omega)\times\omega)$. Define $f:\tilde{V}\rightarrow \tilde{V}$
by $f(x,y)=(y,x)$ for each $x,y\in V$. For each $\gamma\in
C([0,\tau],\Sp(V))$, we define $R(\gamma)\in C([0,\tau],\Sp(V))$ by
$R(\gamma)(t)=\gamma(\tau-t)$ for $t\in[0,\tau]$.

\begin{lemma}\label{l:i-Reverse}
Let $f$ and $R$ be the maps defined above. For each $\gamma\in
C([0,\tau],\Sp(V))$ and $W\in\mathcal{L}(\tilde{V})$, we have
\begin{equation}\label{e:i-Reverse}
    i_W(\gamma(0)^{-1}R(\gamma)\gamma(\tau)^{-1};-\omega)=i_{f(W)}(\gamma).
\end{equation}
\end{lemma}
\begin{proof}
Define $T(\gamma)\in C([0,2\tau],\Sp(V))$ by
\begin{equation*}
    T(\gamma)(t)=\begin{cases}
    \gamma(t) &\text{if $t\in [0,\tau]$},\\
    \gamma(\tau)\gamma(t-\tau)^{-1}\gamma(0)&\text{if $t\in[\tau,2\tau]$}.
    \end{cases}
\end{equation*}
Then $T(\gamma)$ is a contractible loop with a contraction
$T(\gamma(s\cdot))$ for $s\in[0,1]$. So we have
\begin{equation*}
i_{f(W)}(T(\gamma))=i_{f(W)}(\gamma)+i_{f(W)}(\gamma(\tau)\gamma^{-1}\gamma(0))=0.
\end{equation*}
By the definition of Maslov-type index and
$f\in\Sp((\tilde{V},-\tilde{\omega}),(\tilde{V},\tilde{\omega}))$ we
have
\begin{equation*}
    i_{W}(\gamma(0)^{-1}R(\gamma)\gamma(\tau)^{-1};-\omega)=i_{f(W)}(\gamma(\tau)R(\gamma)^{-1}\gamma(0))=-i_{f(W)}(\gamma(\tau)\gamma^{-1}\gamma(0)).
\end{equation*}
So (\ref{e:i-Reverse}) holds.
\end{proof}

Let $\dim V=2n$, $\alpha_1, \alpha_2$ be two Lagrangian subspaces of
$(V,\omega)$ with $V=\alpha_1\oplus\alpha_2$. Set
$\tilde{\alpha}_j:=\alpha_j\times\alpha_j$ for $j=1,2$. Then we have
$\tilde{\alpha}_j\in\mathcal{L}(\tilde{V})$. Set
$N=(-I_{\alpha_1})\oplus I_{\alpha_2}$. Then $N^2=I_V$ and
$\omega(Nx,Ny)=-\omega(x,y)$ for each $x,y\in V$.

\begin{corollary}\label{c:i-N}
Let $f$ and $R$ be the maps defined above. For each $\gamma\in
C([0,\tau],\Sp(V))$ and $W\in\mathcal{L}(\tilde{V})$, we have
\begin{equation}\label{e:i-N}
i_{W}(N\gamma(0)^{-1}R(\gamma)\gamma(\tau)^{-1}N)=i_{f(\diag(N,N)W)}(\gamma).
\end{equation}
\begin{proof}
By Lemma \ref{l:calculate-Maslov-multiply}, Lemma \ref{l:i-Reverse}
and $N\in\Sp((V,\omega),(V,-\omega))$ we have
\begin{equation*}
i_W(N\gamma(0)^{-1}R(\gamma)\gamma(\tau)^{-1}N)=i_{\diag(N,N)W}(\gamma(0)^{-1}R(\gamma)\gamma(\tau)^{-1};-\omega)=i_{f(\diag(N,N)W)}(\gamma).
\end{equation*}
\end{proof}
\end{corollary}

Note that we have
$\Graph(I_V)\in\mathcal{L}_0(\tilde{\alpha}_1,\tilde{\alpha}_2)$,
where $\mathcal{L}_0(\cdot,\cdot)$ is defined by
(\ref{e:definition-L0}). We denote by $p_j$ the projection of
$\tilde{V}$ onto $\tilde{\alpha}_j$ induced by
$\tilde{V}=\tilde{\alpha}_1\oplus\tilde{\alpha}_2$ for $j=1,2$. Then
for each $\mu\in\mathcal{L}_0(\tilde{\alpha}_1,\tilde{\alpha}_2)$,
we have
\begin{equation}\label{e:definition-p-S}
p_j(\mu)=\tilde{\alpha}_j\cap\mu,\quad\Graph(I_V)\cap\mu=S_1(\mu)\oplus
S_2(\mu),
\end{equation}
where we denote by $S_j(\mu):=p_j(\mu)\cap\Graph(I_{\alpha_j})$.
Moreover, the map $p_j$ induces a bijection between
$\mathcal{L}_0(\tilde{\alpha}_1,\tilde{\alpha}_2)$ and the
Grassmannian of $\tilde{\alpha}_j$ for each $j=1,2$.

The following proposition is useful in comparing Maslov-type indices
with a class of Lagrangian boundary conditions.

\begin{proposition}\label{p:s-Lo-p}
Let $\lambda,\mu_1,\mu_2$ be three Lagrangian subspaces of
$\tilde{V}$. Assume that
$\mu_1,\mu_2\in\mathcal{L}_0(\tilde{\alpha}_1,\tilde{\alpha}_2)$ and
$p_2(\mu_1)\subset p_2(\mu_2)$. Then we have
\begin{align}
% \nonumber to remove numbering (before each equation)
  \label{e:s-Lo-p-1}s(\Graph(I_V),\lambda;&\mu_1,\mu_2) = \dim S_2(\mu_2)-\dim S_2(\mu_1)-i(\lambda,\mu_1,\mu_2) \\
\label{e:s-L0-p-2}   &=i(\lambda,\mu_2,\mu_1)-\dim S_1(\mu_1)+\dim S_1(\mu_2) \\
\label{e:s-L0-p-ineq}  &\geq\dim
\frac{\lambda\cap\mu_1}{\lambda\cap\mu_1\cap\mu_2}-\dim
  S_1(\mu_1)+ \dim S_1(\mu_2).
\end{align}
\end{proposition}
\begin{proof}
Since $\mu_1,\mu_2\in
\mathcal{L}_0(\tilde{\alpha}_1,\tilde{\alpha}_2)$ and
$p_2(\mu_1)\subset p_2(\mu_2)$, we have $p_1(\mu_1)\supset
p_1(\mu_2)$. So we have $\mu_1+\mu_2=p_1(\mu_1)+p_2(\mu_2)$ and
$\mu_1\cap\mu_2=p_1(\mu_2)+p_2(\mu_1)$. Together with
(\ref{e:definition-p-S}) we have
\begin{equation*}
    \begin{split}
      \Graph(I_V)\cap(\mu_1+\mu_2) &=S_1(\mu_1)\oplus S_2(\mu_2)  \\
        & =\Graph(I_V)\cap\mu_1+\Graph(I_V)\cap\mu_2,\\
        \Graph(I_V)\cap(\mu_1\cap\mu_2)&=S_1(\mu_2)\oplus
        S_2(\mu_1).
    \end{split}
\end{equation*}
Hence there holds that $\Graph(I_V)\in\mathcal{L}_0(\mu_1,\mu_2)$.
By Corollary \ref{c:kerQ-0}, we have
\begin{equation*}
\begin{split}
i(\Graph(I_V),&\mu_1,\mu_2)=\dim (\Graph(I_V)\cap\mu_2)-\dim
(\Graph(I_V)\cap\mu_1\cap\mu_2)\\
&=(\dim S_1(\mu_2)+\dim S_2(\mu_2))-(\dim S_1(\mu_2)+\dim
S_2(\mu_1))\\
&=\dim S_2(\mu_2)-\dim S_2(\mu_1),
\end{split}
\end{equation*}
and similarly $i(\Graph(I_V),\mu_2,\mu_1)=\dim S_1(\mu_1)-\dim
S_1(\mu_2)$. By Theorem \ref{t:calculate-hormander-index} and Lemma
\ref{l:calculate-triple-index} we have
\begin{equation*}
\begin{split}
  s(\Graph(I_V),\lambda;&\mu_1,\mu_2)=i(\Graph(I_V),\mu_1,\mu_2)-i(\lambda,\mu_1,\mu_2)\\
    & = \dim S_2(\mu_2)-\dim S_2(\mu_1)-i(\lambda,\mu_1,\mu_2)\\
    &=-s(\Graph(I_V),\lambda;\mu_2,\mu_1)\quad (\text{by (\ref{e:s-opposite})})\\
    &= i(\lambda,\mu_2,\mu_1)-i(\Graph(I_V),\mu_2,\mu_1)\\
    &=i(\lambda,\mu_2,\mu_1)-\dim S_1(\mu_1)+\dim S_1(\mu_2)\\
    &\geq\dim\frac{\lambda\cap\mu_1}{\lambda\cap\mu_1\cap\mu_2}-\dim
    S_1(\mu_1)+\dim S_1(\mu_2).
\end{split}
\end{equation*}
\end{proof}

\begin{proposition}\label{p:s-beta-alpha}
Let $\lambda$ be a Lagrangian subspace of $\tilde{V}$ and $\beta$ be
a Lagrangian subspace of $V$. Then we have
\begin{equation}\label{e:s-beta-alpha}
s(\Graph(I_V),\lambda;\tilde{\alpha}_1,\beta\times\alpha_1)\geq \dim
\frac{\lambda\cap\tilde{\alpha}_1}{\lambda\cap((\beta\cap\alpha_1)\times
\alpha_1)}+\dim (\beta\cap\alpha_1)-n.
\end{equation}
\end{proposition}
\begin{proof}
We have
$\Graph(I_V)\in\mathcal{L}_0(\beta\times\alpha_1,\tilde{\alpha}_1)$.
In fact, we have
\begin{equation*}
\Graph(I_V)\cap(\beta\times\alpha_1+\tilde{\alpha}_1)=\Graph(I_V)\cap\tilde{\alpha}_1=\Graph(I_V)\cap(\beta\times\alpha_1)+\Graph(I_V)\cap\tilde{\alpha}_1.
\end{equation*}
By Corollary \ref{c:kerQ-0} we have
\begin{eqnarray*}
i(\Graph(I_V),\beta\times\alpha_1,\tilde{\alpha}_1) & =&\dim
(\Graph(I_V)\cap\tilde{\alpha}_1)\\
         &&-\dim
         (\Graph(I_V)\cap(\beta\times\alpha_1)\cap\tilde{\alpha}_1)\\
         &=&n-\dim(\beta\cap\alpha_1).
\end{eqnarray*}
By (\ref{e:s-opposite}), Theorem \ref{t:calculate-hormander-index}
and Lemma \ref{l:calculate-triple-index} we have
\begin{equation*}
\begin{split}
s(\Graph(I_V),\lambda;&\tilde{\alpha}_1,\beta\times\alpha_1)=i(\lambda,\beta\times\alpha_1,\tilde{\alpha}_1)-i(\Graph(I_V),\beta\times\alpha_1,\tilde{\alpha}_1)\\
&\geq \dim
\frac{\lambda\cap\tilde{\alpha}_1}{\lambda\cap(\beta\times\alpha_1)\cap\tilde{\alpha}_1}-(n-\dim
(\beta\cap\alpha_1))\\
&=\dim\frac{\lambda\cap\tilde{\alpha}_1}{\lambda\cap((\beta\cap\alpha_1)\times\alpha_1)}+\dim
(\beta\cap\alpha_1)-n.
\end{split}
\end{equation*}
\end{proof}

The iteration of a symplectic path $\gamma\in\Pp_{\tau}(V)$ for the
brake symmetry is defined as follows.

\begin{definition}\label{d:brake-iteration} (cf. \cite[(4.3),(4.4)]{LiZh14a}) Given a $\tau>0$, a positive integer $k$, and a path $\gamma\in \Pp_{\tau}(V)$, we define \begin{equation}\label{e:brake-iteration1}
\gamma^{(k)}(t)=\begin{cases}\gamma(t-2j\tau)(\gamma(2\tau))^j,&t\in [2j\tau,(2j+1)\tau],j\in[0,\frac{k-1}{2}],j\in\ZZ,\\
N\gamma(2j\tau-t)N(\gamma(2\tau))^j,&t\in [(2j-1)\tau,2j\tau],j\in[1,\frac{k}{2}],j\in\ZZ,
\end{cases}
\end{equation}
where $\gamma(2\tau)=N\gamma(\tau)^{-1}N\gamma(\tau)$. We call
$\gamma^{(k)}$ the \textit{$k$-th $N$-brake iteration} of $\gamma$.
The map $\gamma(2\tau)$ is called the \textit{Poincar\'{e} map} of
$\gamma$ at $2\tau$.
\end{definition}

\begin{theorem}\label{t:iteration-alpha}
Let $\tilde{\gamma}\in\mathcal{P}_{\tau_1+\tau_2}(V)$ be the
iteration of $\gamma_1\in\mathcal{P}_{\tau_1}(V)$ and
$\gamma_2\in\mathcal{P}_{\tau_2}(V)$. Then we have
\begin{equation}\label{e:iteration-alpha}
i_{\tilde{\alpha}_2}(\tilde{\gamma})-i_{\tilde{\alpha}_2}(\gamma_1)-i_1(\gamma_2)\in[\nu_1(\gamma_2)-2n,0].
\end{equation}
\end{theorem}
\begin{proof}
Set $M_j:=\gamma_j(\tau_j)$ for $j=1,2$. By Definition
\ref{d:Maslov-type-index}, Lemma \ref{l:calculate-Maslov-multiply},
Definition \ref{d:Hormander-index} and Corollary \ref{c:s-2-repeat}
we have
\begin{equation*}
\begin{split}
i_{\tilde{\alpha}_2}(\tilde{\gamma})&-i_{\tilde{\alpha}_2}(\gamma_1)-i_1(\gamma_2)=i_{\tilde{\alpha}_2}(\gamma_2M_1)-i_1(\gamma_2)\\
&=i_{(M_1\alpha_2)\times\alpha_2}(\gamma_2)-i_1(\gamma_2)\\
&=s(\Graph(I_V),\Graph(M_2);\Graph(I_V),(M_1\alpha_2)\times\alpha_2)\\
&=-i(\Graph(M_2),\Graph(I_V),(M_1\alpha_2)\times\alpha_2)\in[\nu_1(\gamma_2)-2n,0].
\end{split}
\end{equation*}
\end{proof}

Then we get the following iteration inequality with brake symmetry
(cf. \cite[Theorem 2.4. $1^{\circ}$]{Li10}).

\begin{corollary}\label{c:i-brake-symmetry}
For each $\gamma\in\mathcal{P}_{\tau}(V)$ and $k\in\NN$, we have
\begin{equation}\label{e:i-brake-symmetry}
\begin{split}
i_{\tilde{\alpha}_2}(\gamma^{(k)})\in &\
i_{\tilde{\alpha}_2}(\gamma^{(k-2[\frac{k-1}{2}])})+[\frac{k-1}{2}]i_1(\gamma^{(2)})\\
&+[[\frac{k-1}{2}](\nu_1(\gamma^{(2)})-2n),\nu_1(\gamma^{(2)})-\nu_1(\gamma^{(2[\frac{k-1}{2}])})].
\end{split}
\end{equation}
\end{corollary}
\begin{proof}
By Corollary \ref{c:i-N} and Lemma \ref{l:calculate-Maslov-multiply}
we have
\begin{equation}\label{e:i-brake-2}
i_1(NR(\gamma)\gamma(\tau)^{-1}N)=i_1(\gamma),
\end{equation}
\begin{equation}\label{e:i-brake-3}
i_1(NR(\gamma)\gamma(\tau)^{-1}N\gamma(\tau))=i_1(\gamma(\tau)NR(\gamma)\gamma(\tau)^{-1}N)=i_1(\gamma
N\gamma(\tau)^{-1}N).
\end{equation}
By (\ref{e:i-brake-2}), (\ref{e:i-brake-3}) and path additivity of
Maslov-type index, we have
\begin{equation}\label{e:i-brake-3-2}
i_1((\gamma^{(3)}|_{[\tau,3\tau]})\gamma(\tau)^{-1})=i_1(\gamma^{(2)}).
\end{equation}
Note that $\gamma^{(2[\frac{k-1}{2}])}$ is the $[\frac{k-1}{2}]$-th
iteration of $\gamma^{(2)}$. By (\ref{e:i-brake-3-2}), Theorem
\ref{t:iteration-alpha} and Theorem
\ref{t:itertion-inequalities-any} we have
\begin{eqnarray*}
i_{\tilde{\alpha}_2}(\gamma^{(k)})&\in&i_{\tilde{\alpha}_2}(\gamma^{(k-2[\frac{k-1}{2}])})+i_1(\gamma^{(2[\frac{k-1}{2}])})+[\nu_1(\gamma^{(2[\frac{k-1}{2}])})-2n,0]\\
&\subset&i_{\tilde{\alpha}_2}(\gamma^{(k-2[\frac{k-1}{2}])})+[\frac{k-1}{2}]i_1(\gamma^{(2)})\\
&&[[\frac{k-1}{2}](\nu_1(\gamma^{(2)})-2n),\nu_1(\gamma^{(2)})-\nu_1(\gamma^{(2[\frac{k-1}{2}])})].
\end{eqnarray*}
\end{proof}

We have the following iteration inequality with focal-type boundary
value condition.

\begin{theorem}\label{t:focal-type-boundary}
Let $\alpha$ be a Lagrangian subspace of $V$. Let
$\tilde{\gamma}\in\mathcal{P}_{\tau_1+\tau_2}(V)$ be the iteration
of $\gamma_1\in\mathcal{P}_{\tau_1}(V)$ and
$\gamma_2\in\mathcal{P}_{\tau_2}(V)$. Set $M_j:=\gamma_j(\tau_j)$
for $j=1,2$. Then we have
\begin{equation}\label{e:focal-type-boundary}
\begin{split}
i_{\alpha\times\alpha_1}(\tilde{\gamma})\, -&\, i_{\alpha\times\alpha_1}(\gamma_1)-i_{\tilde{\alpha}_1}(\gamma_2)\\
\geq& \dim
\frac{\Graph(M_2)\cap\tilde{\alpha}_1}{\Graph(M_2)\cap((\alpha_1\cap(M_1\alpha))\times\alpha_1)}\\
&+ \dim ((M_1\alpha)\cap\alpha_1)-n.
\end{split}
\end{equation}
\end{theorem}
\begin{proof}
By Definition \ref{d:Maslov-type-index}, Lemma
\ref{l:calculate-Maslov-multiply}, Definition
\ref{d:Hormander-index} and Proposition \ref{p:s-beta-alpha} we have
\begin{equation*}
    \begin{split}
    i_{\alpha\times\alpha_1}(\tilde{\gamma})\,-&\,i_{\alpha\times\alpha_1}(\gamma_1)-i_{\tilde{\alpha}_1}(\gamma_2)=i_{\alpha\times\alpha_1}(\gamma_2M_1)-i_{\tilde{\alpha}_1}(\gamma_2)\\
    =&\,i_{(M_1\alpha)\times\alpha_1}(\gamma_2)-i_{\tilde{\alpha}_1}(\gamma_2)\\
    =&\,s(\Graph(I_V),\Graph(M_2);\tilde{\alpha}_1,(M_1\alpha)\times\alpha_1)\\
\geq& \dim
\frac{\Graph(M_2)\cap\tilde{\alpha}_1}{\Graph(M_2)\cap((\alpha_1\cap(M_1\alpha))\times\alpha_1)}\\
&+ \dim ((M_1\alpha)\cap\alpha_1)-n.
    \end{split}
\end{equation*}
\end{proof}

The above theorem has the following two important corollaries.

The following corollary generalizes and strengthens \cite[Theorem
13.5.4]{Lo02}. Note that the inequality
$i_{\tilde{\alpha}_1}(\gamma)\geq n$ holds there.

\begin{corollary}\label{c:i-alpha2-alpha1}
For each $k\in\NN$, we have
\begin{equation}\label{e:i-alpha2-alpha1}
i_{\tilde{\alpha}_2}(\gamma^{(k)})\geq
k(i_{\tilde{\alpha}_1}(\gamma))+\nu_{\tilde{\alpha}_1}(\gamma)-n).
\end{equation}
\end{corollary}
\begin{proof}
Set $M:=\gamma(\tau)$. If $x,Mx\in\alpha_1$, we have $NM^{-1}NMx=x$.
So we have
$\nu_{\tilde{\alpha}_1}(\gamma^{(j)})\geq\nu_{\tilde{\alpha}_1}(\gamma)$
for each $j\in\NN$. Let the map $R$ be defined in Lemma
\ref{l:i-Reverse}. By Corollary \ref{c:i-N} we have
\begin{equation}\label{e:i-alpha1}
i_{\tilde{\alpha}_1}(NR(\gamma)\gamma(\tau)^{-1}N)=i_{\tilde{\alpha}_1}(\gamma).
\end{equation}
By Proposition \ref{p:s-Lo-p} and Theorem
\ref{t:focal-type-boundary}, we have
\begin{equation*}
\begin{split}
i_{\tilde{\alpha}_2}(\gamma^{(k)})&\geq
i_{\tilde{\alpha}_1}(\gamma^{(k)})+\nu_{\tilde{\alpha}_1}(\gamma^{(k)})-n\\
&\geq
ki_{\tilde{\alpha}_1}(\gamma)+\sum_{j=1}^{k}(\nu_{\tilde{\alpha}_1}(\gamma^{(j)})-n)\\
&\geq
k(i_{\tilde{\alpha}_1}(\gamma)+\nu_{\tilde{\alpha}_1}(\gamma)-n).
\end{split}
\end{equation*}
\end{proof}

The following corollary is different from \cite[Theorem 3.7]{Lo93}.
It is enough to prove \cite[Theorem 4.2]{Lo93}. Note that the
inequality $i_{\tilde{\alpha}_1}(\gamma^{(2)})\geq
i_{\tilde{\alpha}_1}(\gamma)\geq n$ holds there.
\begin{corollary}\label{c:i-alpha2-alpha1-odd}
For each $k\in\NN$, we have
\begin{eqnarray}
i_{\alpha_2\times\alpha_1}(\gamma^{(2k+1)})&\geq&
i_{\tilde{\alpha}_1}(\gamma)-n+k(i_{\tilde{\alpha}_1}(\gamma^{(2)})+\nu_{\tilde{\alpha}_1}(\gamma^{(2)})-n),\label{e:i-alpha2-alpha1-odd}
\\
\nu_{\tilde{\alpha}_1}(\gamma^{(2)})&\geq&\nu_{\alpha_1\times\alpha_2}(\gamma).\label{e:nu-gamma-2}
\end{eqnarray}
\end{corollary}
\begin{proof}
Set $M:=\gamma(\tau)$ and $M_2:=NM^{-1}NM$. Then we have
$(NM_2)^2=I_V$ and $M_2=NM_2^{-1}N$. For each $x\in
M^{-1}\alpha_2\cap\alpha_1$, we have $M_2x=-x\in\alpha_1$. So
(\ref{e:nu-gamma-2}) holds. Similarly we get
$\nu_{\tilde{\alpha}_1}(\gamma^{(2k)})\geq
\nu_{\tilde{\alpha}_1}(\gamma^{(2)})$.

By Proposition \ref{p:s-Lo-p} and Theorem
\ref{t:focal-type-boundary}, we have
\begin{equation*}
\begin{split}
i_{\alpha_2\times\alpha_1}(\gamma^{(2k+1)})&\geq
i_{\tilde{\alpha}_1}(\gamma^{(2k+1)})-n\\
&\geq
\sum_{j=1}^k(i_{\tilde{\alpha}_1}(\gamma^{(2)})+\nu_{\tilde{\alpha}_1}(\gamma^{(2j)})-n)+i_{\tilde{\alpha}_1}(\gamma)-n\\
&\geq
k(i_{\tilde{\alpha}_1}(\gamma^{(2)})+\nu_{\tilde{\alpha}_1}(\gamma^{(2)})-n)+i_{\tilde{\alpha}_1}(\gamma)-n.
\end{split}
\end{equation*}
\end{proof}

\section{The almost existence of mean indices}\label{s:mean-indices}

In this section we follow the lines of \cite[\S I.8]{Ek90} and
prove almost existence of mean indices for given
complete autonomous Hamiltonian system on compact symplectic
manifold with symplectic trivial tangent bundle and given autonomous
Hamiltonian system on regular compact energy hypersurface of
symplectic manifold with symplectic trivial tangent bundle. Our main tools will be
Kingman's subadditive ergodic theorem and Theorem
\ref{t:itertion-inequalities-any}.

Firstly we review Kingman's subadditive ergodic theorem.
\begin{definition}\label{d:subadditive-process}
A family $\textbf{x}=(x_{kl};k,l\in\mathbb{N},k<l)$ of random
variables satisfying the following three conditions is called
\textit{subadditive process}.
\begin{description}
 \item[$(S_{1})$] Whenever $j<k<l$, $x_{j,l}\leq x_{j,k}+x_{k,l}.$
 \item[$(S_{2})$] The process $(x_{k+1,l+1})$ has the same joint distributions as the process $(x_{k,l}).$
 \item[$(S_{3})$] The expectation $g_{l}=\textbf{E}(x_{0,l})$ exists
 and satisfies $g_{l}\geq -Al$ for some constant A and all $l>1$.
 \end{description}
\end{definition}

\begin{theorem}[Kingman's subadditive ergodic theorem (cf.\cite{Ki68})]\label{t:Kingman-subadditive}
If $\textbf{x}$ is a subadditive process, then the finite limit
\begin{equation}\label{e:subadditive-limit}
    F=\lim_{l\rightarrow\infty}x_{0l}/l
\end{equation}
exists with probability one and in $L^{1}$, and
\begin{equation}\label{e:subadditive-expection}
    \textbf{E}(F)=\inf_{l\geq1}g_{l}/l.
\end{equation}
\end{theorem}

We need the following preparations that are known to the experts.
For the sake of completeness, we include a proof here.

Let $(V,\omega)$ be a symplectic Hilbert space of dimension $m$, and
$J\in \GL (V)$ be such that $\omega(x,y)=\langle Jx,y \rangle$ for
each $x,y\in V$. The following lemma is used to compare the
Maslov-type indices.

\begin{lemma}\label{l:compare-Maslov-type}
Let $B_0,B_1\in C([0,\tau],\mathcal{B}^{sa}(V))$ be two paths of
self-adjoint operators on $V$, where $\tau>0$. Let $\gamma_s(t)$,
$t\in [0,\tau]$ be the solution of
\begin{equation}\label{e:linear-Hamilton}
  \dot{\gamma_{s}}=-J^{-1}B_s\gamma_{s},\quad\gamma_{s}(0)=I_V,\quad
  t\in [0,\tau],
\end{equation}
for $s=0,1$. Assume that $B_0(t)\leq B_1(t)$ for each $t\in
[0,\tau]$. Then we have
\begin{equation}\label{e:index-compare}
i_1(\gamma_0)\leq i_1(\gamma_1).
\end{equation}
\end{lemma}
\begin{proof}
Set $\Delta B(t):=B_1(t)-B_0(t)$ and $B_s(t):=B_0(t)+s\Delta B(t)$.
Let $\gamma_s(t)$, $0\leq s\leq1$ be the solution of
\begin{equation}\label{e:linear-Hamilton-1}
 \left\{\begin{aligned}
          \dot{\gamma}_s(t)&=-J^{-1}B_s(t)\gamma_s(t),\\
          \gamma_s(0)&=I_V.
         \end{aligned}
  \right.
\end{equation}

We claim that $Q(\Graph(\gamma_s(\tau)),s)$ is semi-positive for
each $s\in[0,1]$.

The idea comes from \cite[Proposition 4.1]{Du76}. We solve the
variational equation
\begin{equation*}
\left\{\begin{aligned}
  \frac{\partial}{\partial t}\frac{\partial \gamma_s(t)}{\partial s}&=-J^{-1}\Delta B(t)\gamma_s(t)-J^{-1}B_s(t)\frac{\partial \gamma_s(t)}{\partial s},\\
  \frac{\partial \gamma_s(0)}{\partial s}&=0.
   \end{aligned}
  \right.
\end{equation*}
Then we get
\begin{equation*}
  \frac{\partial \gamma_s(t)}{\partial s}=-\int_0^t\gamma_s(t)\gamma_s(\tilde{t})^{-1}J^{-1}\Delta B(\tilde{t})\gamma_s(\tilde{t})d\tilde{t}.
\end{equation*}
For each $u,v\in V$,set
\begin{equation*}
    Q_s(\tau)(u,v):=Q(\Graph(\gamma_s(\tau)),s)((u,\gamma_s(\tau)u),(v,\gamma_s(\tau)v)).
\end{equation*}
By Lemma \ref{l:tangential-map}, we have
\begin{equation*}
\begin{split}
Q_s(\tau)(u,v)&=\langle-J\gamma_s(\tau)^{-1}\frac{\partial\gamma_s(\tau)}{\partial
s}u,v\rangle\\
&=\langle (J\int_0^{\tau}\gamma_s(\tilde{t})^{-1}J^{-1}\Delta
B(\tilde{t})\gamma_s(\tilde{t})d\tilde{t})(u),v\rangle\\
&=\int_0^{\tau}\langle \Delta
B(\tilde{t})\gamma_s(\tilde{t})u,\gamma_s(\tilde{t})v\rangle
d\tilde{t}.
\end{split}
\end{equation*}
So we have $Q_s(\tau)\geq0$. By \cite[Proposition 3.2.11]{BoZh14} we
have $i_1(\gamma_s(\tau);s\in[0,1])\geq0$.

By homotopy invariance property and path additivity property of
Maslov-type index, we have
\begin{equation*}
i_1(\gamma_1)=i_1(\gamma_0)+i_1(\gamma_s(\tau);s\in [0,1])\geq
i_1(\gamma_0).
\end{equation*}
\end{proof}

The following lemma gives Maslov-type indices for a kind of
symplectic paths.

\begin{lemma}\label{l:Maslov-type-eJ}
Let $c\in\RR$ and $\tau\geq0$. Set
$J_1:=(-J^2)^{-\frac{1}{2}}J$. Define $\gamma\in
\mathcal{P}_{\tau}(V)$ by $\gamma(t):=e^{cJ_1t}$. Then we have
\begin{equation}\label{e:Maslov-type-eJ}
i_1(\gamma)=mE(\frac{c\tau}{2\pi}).
\end{equation}
\end{lemma}
\begin{proof}
We have $i_1(\gamma)=0$ if $c\tau=0$. Since $J_1^2=-I_V$, $\ker
(\gamma(t)-I_V)\neq 0$ holds if and only if $\frac{ct}{2\pi}\in\ZZ$.
If $\frac{ct}{2\pi}\in\ZZ$, we have $\gamma(t)=I_V$. By Lemma
\ref{l:tangential-map}, $Q(\Graph\circ\gamma,t)$ is positive
(negative) for each $t\in [0,\tau]$ if $c>0$ $(c<0)$ and $\tau>0$.
By Proposition \ref{p:calculate-maslov-index}, we have
\begin{equation*}
i_1(\gamma)=\sum_{\frac{ct}{2\pi}\in
\ZZ,t\in[0,\tau)}m=mE(\frac{c\tau}{2\pi})
\end{equation*}
if $c>0$ and $\tau>0$, and
\begin{equation*}
i_1(\gamma)=-\sum_{\frac{ct}{2\pi}\in\ZZ,t\in(0,\tau]}m=-m[\frac{-c\tau}{2\pi}]=mE(\frac{c\tau}{2\pi})
\end{equation*}
if $c<0$ and $\tau>0$.
\end{proof}

We come to the almost existence of the mean indices.

\begin{ass}\label{a:V-M-B} We make the following assumptions.
\begin{enumerate}
  \item [(i)] Let $(V,\omega)$ be a symplectic Hilbert space of
  dimension $m$, and $J\in\GL(V)$ be such that $\omega(x,y)=\langle
  Jx,y\rangle$ for each $x,y\in V$.
  \item [(ii)] Let $M$ be a compact Hausdorff space with a flow
  $\varphi:\RR\times M\rightarrow M$ and a $\varphi$-invariant
  measure $\mu$ such that $\mu(M)\in(0,+\infty)$.
  \item [(iii)] Let $B\in C(M,\mathcal{B}^{sa}(V))$ be a continuous
  map.
\end{enumerate}
\end{ass}

\begin{definition}\label{d:Maslov-type-M-varphi}
Assume that Assumption \ref{a:V-M-B} holds. For each $\xi\in M$, we
denote by $\gamma_{\xi}$ the fundamental solution of
$\dot{x}=-J^{-1}B(\varphi(t,\xi))x$. For each $\xi\in M$ and
$\tau>0$, the \textit{Maslov-type index} $i_{\tau}(\xi)$ is defined
by
\begin{equation*}
i_{\tau}(\xi):=i_1(\gamma_{\xi}|_{[0,\tau]}).
\end{equation*}
\end{definition}

The main result of this section is the following.

\begin{theorem}\label{t:mean-index}
Assume that Assumption \ref{a:V-M-B} holds. Then there is a Borelian
function $F:M\rightarrow\RR$ such that
\begin{equation}\label{e:F-invariant}
F\circ\varphi(t,\cdot)=F\quad \forall t\geq0,
\end{equation}
\begin{equation}\label{e:mean-index}
\frac{i_{\tau}}{\tau}\rightarrow F\ \text{when $\tau\rightarrow
+\infty$},
\end{equation}
the convergence being $L^1$ and almost everywhere.
\end{theorem}

Before proving the above theorem, we give the following lemma.

\begin{lemma}\label{l:subadditive}
Define the process $\textbf{x}:=(x_{k,l};k,l\in\NN,k<l)$ by
$x_{0,k}=i_k(\xi)$, $k\in\NN$ and
$x_{k+1,l+1}(\xi)=x_{k,l}(\varphi(1,\xi))$, $k,l\in\NN$, $k<l$. Then
the process $\textbf{x}$ is subadditive.
\end{lemma}
\begin{proof}
We claim that for fixed $\tau>0$, $i_{\tau}(\cdot):M\rightarrow\RR$
is a measurable function.

In fact, by Lemma \ref{l:Maslov-index-neighborhood}, for each
$a\in\RR$, the set $A_a:=\{\xi|i_{\tau}(\xi)\geq a\}$ is open and
hence measurable. So $i_{\tau}(\cdot):M\rightarrow \RR$ is a
measurable function.

Now we check that the conditions $(S_{1})$-$(S_{3})$ in Definition
\ref{d:subadditive-process} of subadditive process hold.

$(S_1)$: By the definition of $\textbf{x}$, for $j<k<l$,
$j,k,l\in\NN$, we have
\begin{equation*}
  x_{j,k}(\xi)=x_{j-1,k-1}(\varphi(1,\xi))=\cdots=x_{0,k-j}(\varphi(j,\xi))=i_{k-j}(\varphi(j,\xi)).
\end{equation*}
Similarly, we have
\begin{eqnarray*}
% \nonumber to remove numbering (before each equation)
  x_{k,l}(\xi)& = x_{0,l-k}(\varphi(k,\xi)) = i_{l-k}(\varphi(k,\xi)),\\
  x_{j,l}(\xi)& = x_{0,l-j}(\varphi(j,\xi)) = i_{l-j}(\varphi(j,\xi)).
\end{eqnarray*}
Since $\varphi$ is a flow, we have
\begin{equation*}
 \varphi(k,\xi)=\varphi(k-j,\varphi(j,\xi)).
\end{equation*}
By the above equalities and (\ref{e:iteration-inequality-two}), we
get $(S_1)$.

$(S_2)$: $S_2$ follows from the definition of $\textbf{x}$ and the
measure-preserving property of $\varphi(1,\cdot)$.

$(S_3)$: Since $M$ is compact, $i_{\tau}(\cdot)$ is integrable on
$M$, and there is a constant $c\geq0$ such that
$B(\xi)\geq-c(-J^2)^{\frac{1}{2}}$ for each $\xi\in M$. We then have
the probability space $(M,\frac{\mu}{\mu(M)})$. By Lemma
\ref{l:compare-Maslov-type} and \ref{l:Maslov-type-eJ}, we have
$i_{\tau}(\xi)\geq -m[\frac{c\tau}{2\pi}]$, and
$$g_l=\int_Mx_{0l}(\xi)\frac{d\mu}{\mu(M)}=\int_Mi_{l}(\xi)\frac{d\mu}{\mu(M)}\geq -\frac{clm}{2\pi},$$
for all $l>1$. So $(S_3)$ is verified.

Finally, we get the subaddtive process $\textbf{x}$ on the space
$(M,\frac{\mu}{\mu(M)})$.
\end{proof}

Now we can apply Theorem \ref{t:Kingman-subadditive} to prove
Theorem \ref{t:mean-index}.

\begin{proof}[Proof of Theorem \ref{t:mean-index}]
By Theorem \ref{t:Kingman-subadditive}, $F:=\lim
_{l\rightarrow\infty,l\in\NN}\frac{i_l}{l}$ exists almost everywhere
and in $L^1$, and $F(\varphi(1,\xi))=F(\xi)$, a.e.

For each nonnegative real number $\tau$, by
(\ref{e:iteration-inequality-two}), we have

\begin{equation}\label{e:i-xi-tau}
  -m+i_{[\tau]}(\xi)+i_{\{\tau\}}(\varphi([\tau],\xi))\leq i_{\tau}(\xi)\leq i_{[\tau]}(\xi)+i_{\{\tau\}}(\varphi([\tau],\xi)).
\end{equation}
Since $M$ is compact, there is a constant $c\geq0$ such that
$-c(-J^2)^{\frac{1}{2}}\leq B(\xi)\leq c(-J^2)^{\frac{1}{2}}$ for
each $\xi\in M$. By Lemma \ref{l:compare-Maslov-type} and (\ref{e:Maslov-type-eJ}), for $t\geq0$ we have
\begin{equation}\label{e:i-xi-t}
-m[\frac{ct}{2\pi}]\leq i_t(\xi) \leq mE(\frac{ct}{2\pi}).
\end{equation}
By (\ref{e:i-xi-tau}) and (\ref{e:i-xi-t}) we have
\begin{equation}\label{e:i-tau-t}
  \frac{i_{[\tau]}(\xi)}{[\tau]}\frac{[\tau]}{\tau} - \frac{m[\frac{c}{2\pi}]+m}{\tau}\leq \frac{i_{\tau}(\xi)}{\tau}\leq  \frac{i_{[\tau]}(\xi)}{[\tau]}\frac{[\tau]}{\tau} + \frac{mE(\frac{c}{2\pi})}{\tau}.
\end{equation}
Since $\frac{i_{[\tau]}(\xi)}{[\tau]}$ converges almost everywhere
and in $L^1$ as $\tau\rightarrow +\infty$, by (\ref{e:i-tau-t})
there holds that $\frac{i_{\tau}(\xi)}{\tau}$ converges almost
everywhere and in $L^1$ as $\tau\rightarrow +\infty$.

Moreover, we claim that
\begin{equation}\label{e:F-invariant-1}
F(\xi)=F(\varphi(t,\xi))\quad \forall t\geq0,\xi\in M.
\end{equation}
In fact,
\begin{eqnarray*}
% \nonumber to remove numbering (before each equation)
  F(\xi) &=& \lim_{\tau\rightarrow+\infty} \frac{i_{\tau}(\xi)}{\tau},\\
  F(\varphi(t,\xi)) &=& \lim_{\tau\rightarrow+\infty}
  \frac{i_{\tau}(\varphi(t,\xi))}{\tau}.
\end{eqnarray*}
By (\ref{e:iteration-inequality-two}) we have
\begin{eqnarray*}
-m&\leq i_{\tau+t}(\xi)-i_{\tau}(\xi)-i_t(\varphi(\tau,\xi))&\leq0,\\
-m&\leq i_{\tau+t}(\xi)-i_t(\xi)-i_{\tau}(\varphi(\tau,\xi))&\leq0.
\end{eqnarray*}
So we have
\begin{equation}\label{e:i-t-tau}
-m\leq
i_{\tau}(\varphi(t,\xi))+i_t(\xi)-i_{\tau}(\xi)-i_t(\varphi(\tau,\xi))\leq
m.
\end{equation}
By (\ref{e:i-xi-t}) and (\ref{e:i-t-tau}) we have
\begin{equation}\label{e:i-t-l-1}
-\frac{2m[\frac{ct}{2\pi}]+m}{\tau}\leq
\frac{i_{\tau}(\varphi(t,\xi))}{\tau}-\frac{i_{\tau}(\xi)}{\tau}\leq
\frac{2mE(\frac{ct}{2\pi})+m}{\tau}.
\end{equation}
By letting $\tau\rightarrow+\infty$, we get (\ref{e:F-invariant-1}).
\end{proof}

Let $(M,\omega)$ be a $C^2$ compact symplectic manifold of dimension
$2n$ with $C^2$ boundary.
Let $H:M\rightarrow\RR$ be a $C^2$
function, called Hamiltonian function, which induces a Hamiltonian
vector field $X_{H}: M\rightarrow TM$ defined by
\begin{equation*}
    \omega (X_H(\xi),Y)=-d_{\xi}H(\xi)Y\quad \text{for}\ Y\in T_{\xi}M.
\end{equation*}
Denote by $\varphi(t,\xi)$ the
Hamiltonian flow on $M$ generated by the vector filed $X_H$, that is,
\begin{equation}\label{eq:hamilton}
  \frac{d}{dt}\varphi(t,\xi)=X_{H}(\varphi(t,\xi)),\ \ \  \varphi(0,\xi)=\xi.
\end{equation}

Assume
that $X_{H}$ is complete, i.e., $\varphi(t,\xi)$ is well-defined for
all $(t,\xi)\in\RR\times M$. Denote by $\mu_M:=\frac{\omega^n}{n!}$
the Liouville form of $(M,\omega)$. By Cartan's formula $L_X\omega=X\lrcorner(d\omega)+d(X\lrcorner\omega)$, where $L_X$ denotes the Lie derivative of vector field $X$ and $X\lrcorner$ denotes the interior multiplication with $X$, we have $\varphi(t,\cdot)^*\omega=\omega$ and then $\mu_M$
is $\varphi$-invariant.

Set $V:=\RR^{2n}$ with the standard symplectic form $\omega_0$ on $\RR^{2n}$, then $M\times V$ is a trivial symplectic vector bundle. We set
\begin{equation*}
  \omega_0(x,y)=\langle J_0x,y\rangle \ \text{for all}\ x,y\in\RR^{2n},
\end{equation*}
where $\langle x,y\rangle=y^Tx$, $x,y\in\RR^{2n}$, $y^T$ is the transpose of vector $y$, and
\begin{equation*}
J_0=\left(
                                                                    \begin{array}{cc}
                                                                      0 & -I_n\\
                                                                      I_n & 0 \\
                                                                    \end{array}
                                                                  \right).
\end{equation*}

In \cite{SaZe92}, the Maslov index was defined for non-degenerate periodic solutions of Hamiltonian systems which are contractible loops on a symplectic manifold
$M$ provided that the first Chern class $c_1(TM)$ of the tangent bundle vanishes over
the 2-dimensional homotopy group $\pi_2(M)$.

In our paper, we assume that there is a $C^1$ symplectic trivialization $f:TM\rightarrow M\times V$ of the tangent bundle $TM$ of $M$, that is, $f$ is a $C^1$ bundle isomorphism and $f^*\omega_0=\omega$, such that the following diagram commutes:
\begin{equation*}
\xymatrix{
  TM \ar[rr]^{f} \ar[dr]_{\pi}
                &  &    M\times V \ar[dl]^{p_1}    \\
                & M                }
                \end{equation*}
  where $\pi:TM\rightarrow M$ is the natural projection map, $p_1$ is the projection on the first factor.

Now we use the $C^1$ symplectic trivialization $f$ to define a linear isomorphism $\Phi(\xi)\in \Hom (V,T_{\xi}M)$ for each $\xi\in M$ by $f(\Phi(\xi)x)=(\xi,x)$ for each $x\in V$. Then we have $\Phi(\xi)^*\omega=\omega_0$ for each $\xi\in M$. 
As in \cite{Le03}, we will often commit the usual mild sin of identifying $T_{\xi}M$ with
its image under the canonical injection $X\in T_{\xi}M\mapsto (\xi,X)\in TM$, and will use any of the notations $(\xi,X)$, $X_{\xi}$, and $X$ for a tangent vector in $T_{\xi}M$, depending on how much emphasis we wish to give the point $\xi\in M$.

Consider the linearized flow along the Hamiltonian flow $\varphi(t, \cdot):M\rightarrow M$ of (\ref{eq:hamilton}), the differential with respect to $\xi\in M$, $d_{\xi}\varphi(t,\xi):T_{\xi}M\rightarrow T_{\varphi(t,\xi)}M$,
$t\in\RR$, we can define a map
$\gamma(t,\xi):T_{\xi}M\rightarrow T_{\xi}M$ by
\begin{equation*}
  \gamma(t,\xi)= \Phi(\xi)\Phi(\varphi(t,\xi))^{-1}d_{\xi}\varphi(t,\xi),
\end{equation*}
where $\Phi(\xi)^{-1}:T_{\xi}M\rightarrow V$ is the inverse of the linear isomorphism $\Phi(\xi)$, for any $\xi\in M$,
then we have $\gamma(t,\xi)^*\omega=\omega$. We can use the symplectic trivialization to identity $T_{\xi}M$ with $V$, for each $\xi\in M$. Then for any fixed $\xi\in M$, we get a path $\gamma_0(t,\xi)$, $t\in\RR$ in $\Sp(V_0,\omega_0)$ by
\begin{equation}\label{def:gamma0}
  \gamma_0(t,\xi):=\Phi(\xi)^{-1}\gamma(t,\xi)\Phi(\xi)=\Phi(\varphi(t,\xi))^{-1}d_{\xi}\varphi(t,\xi)\Phi(\xi).
\end{equation}
For any $x\in V$,
\begin{eqnarray*}
% \nonumber to remove numbering (before each equation)
f(\varphi(t,\xi),d_{\xi}\varphi(t,\xi)\Phi(\xi)x)&=&f(\varphi(t,\xi),\Phi(\varphi(t,\xi))\gamma_0(t,\xi)x) \\
&=&(\varphi(t,\xi),\gamma_0(t,\xi)x).
\end{eqnarray*}

 For any fixed $\xi_0\in M$, choose a chart $(U_0,\psi)$ for $M$ with coordinate functions $(\xi^1,...,\xi^{2n})$, such that $\xi_0\in U_0$. Define a map $\tilde{\psi}:\pi^{-1}(U_0)\rightarrow \psi(U_0)\times \RR^{2n}$ by
 \begin{equation*}
   \tilde{\psi}(\sum_{k=1}^{2n}v^k \frac{\partial}{\partial\xi^k}|_{\xi})=(\xi^1(\xi),...,\xi^{2n}(\xi),v^1,...,v^{2n}),
 \end{equation*}
 where $\frac{\partial}{\partial\xi^k}|_{\xi}\in T_{\xi}M$ are the coordinate vectors associated with the given chart and $(\frac{\partial}{\partial\xi^1},...,\frac{\partial}{\partial\xi^{2n}})$ form a $C^1$ local frame for $TM$ over $U_0$. The $C^1$ map $\tilde{\psi}$ is a bijection onto its image, obviously linear on fibers and satisfies $p_1\circ \tilde{\psi}=\pi$.
 The coordinates $(\xi^1,...,\xi^{2n},v^1,...,v^{2n})$ are called standard coordinates for $TM$.
 Using the local coordinate representation,
the Cauchy problem of the Hamiltonian system (\ref{eq:hamilton}) reduces to an initial value problem of a first order $C^1$ system of ordinary differential equations.
Then the solution $\varphi(t,\xi)$ of (\ref{eq:hamilton}) is $C^1$ on
\begin{equation}\label{eq:phismallt}
  G_0:=\{(t,\xi);t\in (-\delta_{\xi_0},\delta_{\xi_0}),\xi\in U_{\xi_0}\}
\end{equation}
for some positive number $\delta_{\xi_0}$ and open set $U_{\xi_{0}}\subseteq U_0$ such that $\varphi(G_0)\subset U_0$. Denote by $g(t,\xi):=(\varphi(t,\xi),d_{\xi}\varphi(t,\xi)\Phi(\xi)x)\in TM$ for $t\in\RR$, $\xi\in M$, $x\in V$. Then for fixed $\xi\in M$, $g(t,\xi)$ is a curve in $C^1$ manifold $TM$.
We claim that $\frac{d}{dt}g(t,\xi)$ exists and is continuous on $G_0$. Thus for each $\xi\in U_{\xi_0}$, $\frac{d}{dt}g(t,\xi)$ is the velocity vector field of the curve $g(t,\xi), t\in (-\delta_{\xi_0},\delta_{\xi_0})$.
Locally, using the standard coordinates chart $(\pi^{-1}(U_0),\tilde{\psi})$ for $TM$,
\begin{eqnarray*}
% \nonumber to remove numbering (before each equation)
  \tilde{\psi}(\varphi(t,\xi),d_{\xi}\varphi(t,\xi)\Phi(\xi)x) &=&\tilde{\psi}\circ d_{\xi}\varphi(t,\xi)\circ\tilde{\psi}^{-1}\circ\tilde{\psi}(\Phi(\xi)x) \\
   &=& (\psi(\varphi(t,\xi)),(\frac{\partial \varphi^j}{\partial\xi^k}(t,\xi))_{j,k=1}^{2n}\tilde{\psi}(\Phi(\xi)x)),
\end{eqnarray*}
where the matrix $(\frac{\partial \varphi^j}{\partial\xi^k})_{j,k=1}^{2n}$ of partial derivatives is the \textit{Jacobian matrix} of $\psi\circ\varphi\circ\psi^{-1}$.
In fact, although $\varphi$ may not be $C^2$, by the standard coordinates $(\xi^1,...,\xi^{2n},v^1,...,v^{2n})$ representation and the $C^1$ system (\ref{eq:hamilton}),
$(\frac{\partial \dot{\varphi}^j}{\partial\xi^k}(t,\xi))_{j,k=1}^{2n}$ is continuous on
$G_0$ and $\dot{\varphi}(t,\xi)$, $d_{\xi}\varphi(t,\xi)$ exist. Then we have $\frac{d}{dt}(\frac{\partial \varphi^j}{\partial\xi^k}(t,\xi))_{j,k=1}^{2n}$
exists and $\frac{d}{dt}(\frac{\partial \varphi^j}{\partial\xi^k}(t,\xi))_{j,k=1}^{2n}=(\frac{\partial \dot{\varphi}^j}{\partial\xi^k}(t,\xi))_{j,k=1}^{2n}$, since every term in the matrices is just the second-order mixed partial derivative.
So
\begin{eqnarray*}
% \nonumber to remove numbering (before each equation)
  \frac{d}{dt}((\frac{\partial \varphi^j}{\partial\xi^k}(t,\xi))_{j,k=1}^{2n}\tilde{\psi}(\Phi(\xi)x))&=&\frac{d}{dt}(\frac{\partial \varphi^j}{\partial\xi^k}(t,\xi))_{j,k=1}^{2n}\tilde{\psi}(\Phi(\xi)x)  \\
   &=&(\frac{\partial \dot{\varphi}^j}{\partial\xi^k}(t,\xi))_{j,k=1}^{2n}\tilde{\psi}(\Phi(\xi)x),
\end{eqnarray*}
thus $\frac{d}{dt}g(t,\xi)$ exists and is continuous on $G_0$.
 Since $f$ is $C^1$, using the chain rule for total derivatives to the composite curve, we obtain that $\dot{\gamma}_0(t,\xi)$ is continuous on $G_0$.

Since $M$ is compact, $M$ can be covered by a finite number of such neighborhood $U_{\xi_0}$. Let $\delta_0$ denote the smallest of the corresponding positive number $\delta_{\xi_0}$.

For any $\xi\in M$, since $\varphi(0,\xi)=\xi$, we have $\gamma_0(0,\xi)=I_V$.

For any $\xi\in M$, $s,t\geq 0$, since
\begin{equation}\label{eq:varphiflow}
  \varphi(s+t,\xi)=\varphi(s,\varphi(t,\xi)),
\end{equation}
and by the chain rule for total derivatives, we have
\begin{equation}\label{eq:dphitime}
  d_{\xi}\varphi(s+t,\xi)=d_{\zeta}\varphi(s,\zeta)|_{\zeta=\varphi(t,\xi)}d_{\xi}\varphi(t,\xi).
\end{equation}
By (\ref{def:gamma0}) and (\ref{eq:dphitime}), we have
\begin{equation}\label{eq:gammaflow}
  \gamma_0(s+t,\xi)=\gamma_0(s,\varphi(t,\xi))\gamma_0(t,\xi).
\end{equation}

We claim that $\varphi(t,\xi)$ is $C^1$ on $\RR\times M$.
In fact, for any $t\in \RR$, we can choose an $N_0\in\NN$, such that
$\frac{|t|}{N_0}<\delta_0$, then by (\ref{eq:varphiflow}) and (\ref{eq:phismallt}),
we have this claim.

So by (\ref{eq:gammaflow}), $\dot{\gamma_0}(t,\xi)$ exists and is continuous on $\RR\times M$. Since
\begin{equation*}
  \gamma_0(t,\xi)^{T}J_0\gamma_0(t,\xi)=J_0,\ \text{for}\ t\in\RR,\xi\in M,
\end{equation*}
simple computations show that $-J_0\dot{\gamma_0}(t,\xi)\gamma_0(t,\xi)^{-1}$ are symmetric matrices. Define $B:M\rightarrow \mathcal{B}^{sa}(V)$ by  $B(\xi)=-J_0\dot{\gamma_0}(t,\xi)|_{t=0}$, then we get a $B\in C(M,\mathcal{B}^{sa}(V))$. Altogether,
we have for fixed $\xi$, $\gamma_0(t,\xi)$ is the fundamental solution of
$\dot{x}=-J_0^{-1}B(\varphi(t,\xi))x$. By Theorem \ref{t:mean-index} we have
\begin{corollary}\label{e:i-manifold-boundary}
Under the above assumptions, the results of Theorem
\ref{t:mean-index} hold.
\end{corollary}

Let $(M,\omega)$ be a $C^3$ compact symplectic manifold of dimension
$2n$. Set $V:=\RR^{2n}$. Assume that there is a $C^1$ symplectic
trivialization $f:TM\rightarrow M\times V$ of the tangent bundle
$TM$ of $M$. Let $\Sigma$ be an orientable compact $C^3$
hypersurface of $M$. By \cite[\S4.2, p. 114]{HoZe94}, there exists a
function $H\in C^2(U,\RR)$ defined on an open neighborhood $U$ of
$\Sigma$ representing $\Sigma=\{x\in U; H(x)=0\}$ and satisfying
$dH\neq0$. Moreover, there exists an $\epsilon>0$ and a
diffeomorphism
\begin{equation*}
    \psi:\Sigma\times I\rightarrow U\subset M
\end{equation*}
such $U=H^{-1}(I)$, and
$$\psi(x,0)=x,\ H(\psi(x,t))=t,\ \text{for} \ (x,t)\in\Sigma\times I, $$
where $I=(-\epsilon,\epsilon)$.

%Define $B:M\rightarrow \mathcal{B}^{sa}(V)$ by $\langle
%B(\xi)u,v\rangle=(d^2H)(\xi)(f^{-1}(\xi,u),f^{-1}(\xi,v))$ for each
%$\xi\in M$ and $u,v\in V$.
$H$ induces a Hamiltonian vector field
$X_H$ on $U$. Denote by $\varphi(t,\xi)$ the Hamiltonian flow
generated by $X_H$, then $X_H$ is complete on $U$, and $\Sigma$ is
an invariant subset of the flow $\varphi$. Define $B:M\rightarrow \mathcal{B}^{sa}(V)$ as in Corollary \ref{e:i-manifold-boundary}.
 Denote by
$\mu_M:=\frac{\omega^n}{n!}$ the Liouville form of $(M,\omega)$.

\begin{lemma}\label{l:invariant-measure}
Under the above assumptions, there exists a $(2n-1)$-form $\mu_{\Sigma}$
on $\Sigma$ such that the form $\mu_{\Sigma}$ is
$\varphi$-invariant.
\end{lemma}
\begin{proof}
Let $x=(x_1,x_2,...,x_{2n})\in U'\subset R^{2n}$ be a local
coordinate of $\xi\in U$. Then there exists a positive function $a$ on
$U'$ such that $\mu_M=adx_1\wedge dx_2...\wedge dx_{2n}$.

Let
$\mu=\sum_{i=1}^{2n}(-1)^{i-1}a_{i}dx_{1}\wedge...\wedge\widehat{dx_{i}}...\wedge
dx_{2n}$ be a $(2n-1)$-form defined in the chart such that
\begin{equation}\label{e:measure}
\mu_M=dH\wedge\mu,
\end{equation}
where the hat indicates omitted argument. Denote by
$H_{x_j}:=\frac{\partial H}{\partial x_j}$. Then we have
$$\sum_{j=1}^{2n}H_{x_{j}}a_{j}=a.$$

Denote by $\iota : \Sigma \hookrightarrow M$ the inclusion map. We claim
that, if there are two $(2n-1)$-forms $\mu_{1}$, $\mu_{2}$ which
satisfy (\ref{e:measure}), then $\iota^{*}\mu_{1}$ =
$\iota^{*}\mu_{2}$. In fact, $H(x)=0$ holds for $x\in\Sigma$, so we
have $\sum_{j=1}^{2n}H_{x_{j}}dx_{j}=0$. Since $H'(x)\neq0$, without
generality, we assume $H_{x_{1}}\neq0$, then we have
$\iota^{*}\mu=\frac{a}{H_{x_{1}}}dx_{2}\wedge...\wedge
dx_{2n}|_{\Sigma}$, which is unique determined by $H$. The argument
also shows the existence of $\mu$ in chart.

By patching the local defined $\iota^*\mu$, we obtain a globally
defined $(2n-1)$-form $\mu_{\Sigma}:=\iota^*\mu$ on $\Sigma$.

We have $\varphi(t,\cdot)^*\mu_M=\mu_M$ by Liouville theorem and
$\varphi(t,\cdot)^*dH=dH$ since $H\circ\varphi(t,\cdot)=H$. So we
have
$\varphi(t,\cdot)^*\mu_M=(\varphi(t,\cdot)^*dH)\wedge(\varphi(t,\cdot)^*\mu)$,
and
\begin{equation}\label{e:mrasure-1}
\mu_M=dH\wedge(\varphi(t,\cdot)^*\mu).
\end{equation}
By the uniqueness of $\iota^*\mu$ we have
$\iota^{*}\varphi(t,\cdot)^*\mu=\iota^{*}\mu$. Since
$\varphi(t,\cdot)\circ \iota=\iota \circ \varphi(t,\cdot)$ on
$\Sigma$, we finally get
$$\varphi(t,\cdot)^*\iota^*\mu=\iota^*\mu.$$
\end{proof}

By Theorem \ref{t:mean-index} we have
\begin{corollary}\label{c:i-compact-manifild}
Under the above assumptions of Lemma \ref{l:invariant-measure},
there is a Borelian  function $F:\Sigma\rightarrow\RR$ such that
\begin{equation}\label{e:F-invariant-2}
F\circ\varphi(t,\cdot)=F\quad \forall t\geq0,
\end{equation}
\begin{equation}\label{e:i-mean-index}
\frac{i_{\tau}}{\tau}\rightarrow F\ \text{when $\tau\rightarrow
+\infty$},
\end{equation}
the convergence being $L^1$ and almost everywhere. Here the flow
$\varphi$ and the function $i_{\tau}$ are defined on $\Sigma$.
\end{corollary}

\noindent\bf{\footnotesize Acknowledgements}\quad\rm
{\footnotesize We would like to thank the referees of this paper for their critical reading and very helpful comments and suggestions. The authors were partially supported by NSFC (No.11221091 and No.11471169), LPMC of MOE of China.}\\[4mm]

\noindent{\bbb{References}}

\end{document}